\documentclass[11pt,a4paper]{article}
\usepackage[margin=25mm,top=22mm,bottom=23mm]{geometry}
\usepackage[T1]{fontenc}
\usepackage{lmodern,microtype,amsmath,amssymb,amsthm}
\usepackage[colorlinks=true,linkcolor=blue!45!black,urlcolor=blue!45!black,citecolor=blue!45!black]{hyperref}
\usepackage{xcolor}
\newtheorem{theorem}{Theorem}
\newtheorem{proposition}[theorem]{Proposition}
\newtheorem{lemma}{Lemma}[section]
\newtheorem{claim}[lemma]{Claim}
\newtheorem*{claim*}{Claim}
\theoremstyle{remark}
\newtheorem{remark}[lemma]{Remark}
\newcommand{\R}{\mathbb R}
\newcommand{\E}{\mathcal E}
\newcommand{\X}{X^1(\R)}
\newcommand{\dd}{\,dx}
\allowdisplaybreaks[2]
\title{\vspace{-2.2em}Global well-posedness for the focusing Gross-Pitaevskii equation in the one-dimensional Zhidkov space}
\author{Agissilaos G. Athanassoulis\\[0.3em]
\small University of Dundee\\
\small\texttt{aathanassoulis@dundee.ac.uk}}
\date{}
\hypersetup{pdftitle={Global well-posedness for the focusing Gross\textendash Pitaevskii equation in the one-dimensional Zhidkov space},pdfauthor={Agissilaos Athanassoulis}}
\begin{document}
\maketitle\vspace{-1.1em}

\begin{abstract}
We show global well-posedness for the focusing Gross-Pitaevskii equation in one spatial dimension, with initial data in the Zhidkov space, with no smallness condition. For the usual nonlinear Schr\"odinger equation with decaying initial data, a mass-energy argument is used to preclude finite-time blow-up. Here both mass and total energy are infinite in general, but we put together a renormalized version of that bound by using the net mass transfer and the kinetic energy, which we show are a priori bounded. Global well-posedness on the energy space and on $1+H^1$ follows. The argument extends to other problems with smooth real density responses, including the competing cubic-quintic equation, and exponents in the quadratic-cubic range.

\medskip\noindent\textbf{MSC 2020:} 35Q55 (primary), 35A01, 35B45, 35B30 (secondary).
\end{abstract}

\section{Introduction}

The focusing Gross--Pitaevskii equation on the real line,
\begin{equation}\label{eq:gp}
 iu_t+u_{xx}+q(|u|^2-1)u=0,\qquad u(0)=u_0,\qquad q>0,
\end{equation}
is often posed with initial data that do not decay at infinity. Its name originates in the work on Bose--Einstein condensates \cite{gross,pitaevskii}. Another classical problem linked with the same equation is the modulation instability of ocean waves \cite{BenjaminFeir,Zakharov}. In both cases, the fact that initial data do not decay at infinity was central.

To describe the initial data, we use the functional frameworks for non-decaying data developed by Zhidkov, Gallo and G\'erard \cite{gallo,gerard,zhidkov}, working in the three nested spaces
\[
 1+H^1(\R)\subsetneq \E_1\subsetneq X^1(\R).
\]
The affine space $1+H^1$ consists of $H^1$ perturbations of the constant field $1$ and carries the distance $d_{1+H^1}(v,w)=\|v-w\|_{H^1}$. By phase invariance, any fixed constant of unit modulus can be used in its place.
On the line we equip G\'erard's energy space \cite{gerard} with the metric $d_{\E_1}$,
\begin{equation}\label{eq:energy-space}
 \begin{aligned}
 \E_1 &= \{v\in H^1_{\mathrm{loc}}(\R;\mathbb C): v_x\in L^2(\R),\ |v|^2-1\in L^2(\R)\},\\
 d_{\E_1}(v,w) &= \|v-w\|_\infty+\|v_x-w_x\|_2+\||v|^2-|w|^2\|_2 .
 \end{aligned}
\end{equation}
Functions in the energy space may have different phases at the two ends of the line, or no limiting phase at all \cite{bgss-black}.
The first-order Zhidkov space \cite{gallo,zhidkov} requires only boundedness and a square-integrable derivative,
\begin{equation}\label{eq:zhidkov-space}
 \begin{aligned}
 X^1(\R) &= \{v\in L^\infty(\R;\mathbb C):v_x\in L^2(\R)\},\\
 d_{X^1}(v,w) &= \|v-w\|_\infty+\|v_x-w_x\|_2 .
 \end{aligned}
\end{equation}

Theorem~\ref{thm:main} covers global well-posedness on each of the three spaces, with explicit norm-growth bounds. We use a novel net mass transfer argument to preclude finite-time blow-up. Section~\ref{sec:applications} extends the argument to smooth density responses, including the competing cubic--quintic equation, and to focusing powers $2\le p<3$.

\subsection{State of the art}\label{sec:related}
Zhidkov developed the bounded, finite-gradient framework on the line \cite{zhidkov}. Gallo established the Schr\"odinger group and local theory on Zhidkov spaces \cite{gallo}, and G\'erard developed the energy metric and global theory for the defocusing version of equation \eqref{eq:gp} ($q<0$) in dimensions two and three \cite{gerard}. We use these results heavily.

For the focusing sign, local theories include de Laire's theorem on $\phi+H^1$ \cite{delaire} and Mu\~noz's $H^s$ perturbations of the Stokes wave for $s>1/2$ \cite{munoz}. Hennig et al.\ establish local $H^1$ theory for focusing power and saturable laws with nonzero backgrounds \cite{hennig-proximity}, followed by quantitative lifespan estimates \cite{hennig}. 

For the quadratic equation, Chaichenets et al.\ prove global well-posedness for both signs in $H^1(\mathbb T)+L^2(\R)$, with an exponential perturbation-mass bound \cite{quadratictooth}. For the cubic equation, inverse scattering treats localized perturbations of nonzero endstates \cite{bilman,biondini,biondini-mantzavinos}. Recently, Kotani et al.\ use the integrable NLS hierarchy to obtain global well-posedness for both signs in $M^5_{\infty,1}(\R)$, which includes $1+H^6$ data and smooth periodic fields, with implicit bounds on finite time intervals in the preprint \cite{kotani}. Initial data in $X^1(\mathbb{R})$ need not be in $M^5_{\infty,1}(\R),$ and vice versa: 
the periodic backgrounds allowed in $M^5_{\infty,1}(\R)$ lie outside $X^1(\R)$ in general.

\subsection{Main result and the heuristic behind it}\label{sec:main}
\begin{theorem}[Focusing Gross--Pitaevskii on the line]\label{thm:main}
For every $q>0$ and $u_0\in\X$, \eqref{eq:gp} has a unique global mild solution $u\in C(\R;\X)$. For all $t\in\R$, the solution satisfies the bounds
\begin{align}
 \|u_x(t)\|_2 &\le \|\partial_x u_0\|_2\bigl(2e^{q\|u_0\|_\infty^2|t|}-1\bigr),\label{eq:cubic-kinetic-growth}\\
 \|u(t)\|_\infty &\le \|u_0\|_\infty
 + C q^{-1/6}\|\partial_x u_0\|_2^{2/3} e^{\frac23 q\|u_0\|_\infty^2|t|},\label{eq:cubic-amplitude-growth}
\end{align}
with a universal $C$. Define the solution map by
\[
 \Phi_t(u_0)=u(t).
\]
It restricts to global flows on $\E_1$ and $1+H^1(\R)$. For $u_0\in\E_1$, the density defect satisfies
\begin{equation}\label{eq:cubic-energy-density}
 \big\||u(t)|^2-1\big\|_2
 \le \left(\big\||u_0|^2-1\big\|_2
 +\sqrt{8/q}\,\|\partial_xu_0\|_2\right)
 e^{q\|u_0\|_\infty^2|t|}.
\end{equation}
For $u_0\in1+H^1(\R)$, the full perturbation norm satisfies
\begin{equation}\label{eq:cubic-affine-growth}
 \begin{aligned}
 \|u(t)-1\|_{H^1}
 &\le \|u_0-1\|_{H^1}\\
 &\quad+\left[\big\||u_0|^2-1\big\|_2
 +\bigl(2+\sqrt{8/q}\bigr)\|\partial_xu_0\|_2\right]
 \bigl(e^{q\|u_0\|_\infty^2|t|}-1\bigr).
 \end{aligned}
\end{equation}
On $\E_1,$ the flow conserves the signed energy
$E_q(u)=\|u_x\|_2^2-(q/2)\||u|^2-1\|_2^2$.

Moreover, we have continuous dependence on the initial data in bounded sets for finite time, for each of the three spaces. More specifically, let $Y$ be one of $1+H^1,$ $\E_1,$ $X^1(\R);$  let $R>0$, and denote $B_R=\{v\in Y:d_Y(v,1)<R\}$. Then, for every $u_0\in B_R$,
\begin{equation}\label{eq:cubic-metric-growth}
 d_Y(\Phi_t u_0,1)\le C_{q,R}\,e^{q\|u_0\|_\infty^2 |t|},
\end{equation}
and for any $u_0,v_0\in B_R$ and $T>0$,
\begin{equation}\label{eq:cubic-metric-dependence}
 \sup_{|t|\le T}d_Y(\Phi_t u_0,\Phi_t v_0)\le C_{q,R,T}\,d_Y(u_0,v_0).
\end{equation}
\end{theorem}

\begin{remark}[Idea of the proof: net mass transfer is finite]\label{rem:kinetic-energy}
For the one-dimensional focusing cubic NLS with $H^1$ data, conserved mass and energy combine to preclude blow-up. Here the total mass is infinite; we keep track of the net mass transfer instead.  
Let us start with Gallo's local solution in $X^1$ on a maximal time interval $I\ni 0$ \cite{gallo}.

Denote $n=|u|^2$, $n_0=|u_0|^2$ and $j=2\operatorname{Im}(\bar u u_x)$. The continuity equation gives
\[
 F(t)=-\int_0^t j(s)\,ds,\qquad F_t=-j,\qquad F_x=n-n_0.
\]
The net mass gained by $[a,b]$ is $F(t,b)-F(t,a)$; we call $F$ the \emph{net mass transfer function}. All quantities and integrated identities below are justified for mild $X^1$ solutions throughout their local lifespan, before any hypothetical blow-up. The details are in the proof, in Section~\ref{sec:cubic}.

Testing $F_t=-j$ against $F$ gives
\begin{equation}\label{eq;starstar}
\left|\frac{d}{dt}\|F(t)\|_2^2\right|\le 4\|Fu\|_2\,\|u_x(t)\|_2.
\end{equation}
The {\em mass-weighted cancellation}
\begin{equation}\label{eq:mass-weighted-cancellation}
\int(n-n_0)F^2=\int(F^3)_x/3=0
\implies \|Fu\|_2 = \|Fu_0\|_2
\end{equation}
yields the bound $\|Fu\|_2\le \|u_0\|_\infty \|F(t)\|_2$. This controls $\|Fu\|_2$ in \eqref{eq;starstar}, but leaves the kinetic norm $\|u_x(t)\|_2$ uncontrolled. We therefore replace $u_x$ in the current by the modified derivative
\[
 Y=u_x+\tfrac q2Fu,\qquad j=2\operatorname{Im}(\bar uY).
\]
Note that any real-valued function $g$ would yield
a $\tilde{Y}=u_x + gu$ with $j=2\operatorname{Im}(\bar u\tilde{Y}).$ The reason for our particular choice of $Y$
is the identity \eqref{eq:cubic-square},
\[
 \|Y\|_2=\|\partial_xu_0-\tfrac q2Fu_0\|_2 \le \|\partial_xu_0\|_2+\tfrac q2\|u_0\|_\infty \|F(t)\|_2.
\]
This originates in a computation of the kinetic energy increment $\|u_x(t)\|^2_2 - \|\partial_x u_0\|_2^2$ which can close using $F$ even if the potential energy is infinite.

Now using $Y$ in the current estimate now gives
\[
 \frac{d}{dt}\|F(t)\|_2\le 2\|u_0\|_\infty \|\partial_xu_0\|_2+q\|u_0\|_\infty^2\|F(t)\|_2, \qquad F(0)=0.
\]
Now Gronwall precludes finite-time blow-up for $F$ in $L^2(\mathbb{R});$ more specifically it yields
\begin{equation}\label{eq:star}
\|F(t)\|_2\le \frac{2\|\partial_xu_0\|_2}{q\|u_0\|_\infty}(e^{q\|u_0\|_\infty^2|t|}-1).
\end{equation} 

For the kinetic energy, observe that $u_x=Y-\tfrac q2Fu$. The triangle inequality and the preceding bounds give
\[
 \|u_x(t)\|_2\le \|Y(t)\|_2+\tfrac q2\|F(t)u(t)\|_2
\le \|\partial_xu_0\|_2+q\|u_0\|_\infty \|F(t)\|_2.
\]
Combining this with \eqref{eq:star} we obtain the a priori bound
\[
 \|u_x(t)\|_2\le \|\partial_xu_0\|_2\bigl(2e^{q\|u_0\|_\infty^2|t|}-1\bigr).
\]
The proof then proceeds to control the amplitude $\|u(t)\|_\infty$ and rule out a finite endpoint for the maximal time interval $I$.

\end{remark}

\subsection{Application to other nonlinearities}\label{sec:applications}
The net mass transfer argument can be generalized  to other density responses without requiring finite relative energy. 
Consider a real density response $f;$ the nonlinear Schr\"odinger equation
\begin{equation}\label{eq:nls-unrotated}
 i\psi_t+\psi_{xx}+f(|\psi|^2)\psi=0,\qquad \psi(0)=u_0,
\end{equation}
is equivalent, under the phase rotation $u(t)=e^{-if(1)t}\psi(t)$, to
\begin{equation}\label{eq:nls}
 iu_t+u_{xx}+\bigl(f(|u|^2)-f(1)\bigr)u=0,\qquad u(0)=u_0.
\end{equation}
We write $\mathcal N(z)=(f(|z|^2)-f(1))z$ for the nonlinearity in \eqref{eq:nls}. This normalization makes the unit background stationary. See Lemma~\ref{lem:phase-rotation} for the corresponding equivalence of mild solutions. In particular, the phase rotation (i.e., passing from \eqref{eq:nls-unrotated} to \eqref{eq:nls}) preserves $H^1$ norms of the inhomogeneity and maximal lifespans.

\begin{proposition}[Smooth density responses with bounded slope]\label{prop:bounded-slope}
Let $f\in C^2([0,\infty);\R)$ satisfy
\[
 L=\sup_{s\ge0}|f'(s)|<\infty.
\]
For every $u_0\in\X$, \eqref{eq:nls} has a unique global mild solution $u\in C(\R;\X)$. This flow restricts to globally well-posed flows in $\E_1$ and in $1+H^1(\R)$. The flows in all three spaces are Lipschitz on bounded sets on compact time intervals, in their respective metrics.
%The kinetic norm satisfies the estimate \eqref{eq:main-L} with this value of $L$.
\end{proposition}

The bounded-slope assumption allows unbounded density responses, including $f(s)=qs$, which recovers the cubic equation. In fact, this can be relaxed further:

\begin{proposition}[A one-sided bound on the anti-derivative]\label{prop:anti-derivative-growth}
Let $f\in C^2([0,\infty);\R)$ and suppose that an anti-derivative $V$ of $f$ satisfies
\[
 V(s)\le C_V(1+s^2)\qquad (s\ge0)
\]
for some $C_V>0$. Then \eqref{eq:nls} is globally well-posed in $\X$, $\E_1$ and $1+H^1(\R)$. In each space the flow is Lipschitz on bounded sets on compact time intervals.
In particular, these conclusions hold for the competing cubic--quintic equation
\[
 iu_t+u_{xx}+\bigl[q(|u|^2-1)-h(|u|^4-1)\bigr]u=0,\qquad q,h>0 .
\]
\end{proposition}
Indeed, for the competing cubic-quintic $f(s)=qs-hs^2$ has unbounded slope, while its anti-derivative $V(s)=qs^2/2-hs^3/3$ is bounded above on $[0,\infty)$.
Proposition~\ref{prop:anti-derivative-growth} gives global well-posedness for the cubic--quintic equation for every $q,h>0$ and every $u_0\in X^1(\R)$. This extends Gallo's global existence result \cite{gallo-defocusing}, which applies for initial data in $\E_1$ and $h>q/2$.

For the focusing power $f(s)=qs^{(p-1)/2}$ with $2\le p<3$,  $f'$ is singular at zero, so neither Proposition~\ref{prop:bounded-slope} nor Proposition~\ref{prop:anti-derivative-growth} apply. 

\begin{proposition}[Weaker focusing powers]\label{prop:subcubic-powers}
Let $q>0$ and $2\le p<3$. The equation
\begin{equation}\label{eq:subcubic-power}
 iu_t+u_{xx}+q(|u|^{p-1}-1)u=0,\qquad u(0)=u_0,
\end{equation}
is globally well-posed in $\X$, $\E_1$ and $1+H^1(\R)$. In each space the flow is Lipschitz on bounded sets on compact time intervals. For every $u_0\in\X$ and $t\in\R$,
\begin{equation}\label{eq:subcubic-kinetic}
 \|u_x(t)\|_2\le\|\partial_xu_0\|_2
 \left(2\exp\left(\frac{2q(p-1)}{p+1}\|u_0\|_\infty^{p-1}|t|\right)-1\right).
\end{equation}
\end{proposition}

Section~\ref{sec:lip} proves these propositions using the generalized net mass transfer estimates and relative kinetic balance from Section~\ref{sec:current}.

\paragraph{AI statement.} ChatGPT Astra 6 and 6.2 was used to proofread the manuscript. The author has sole responsibility for the content and correctness of the results.

\section{Functional framework}\label{sec:framework}

Here we collect some standard facts about the spaces we work with. 
All fields are complex valued and all spatial derivatives are understood in weak sense. Throughout the paper the spatial domain is $\R$. We use the spaces and metrics in \eqref{eq:energy-space} and \eqref{eq:zhidkov-space}, with
\[
 \|v\|_{X^1}=\|v\|_\infty+\|v_x\|_2.
\]
The bound $|\partial_x|v||\le |v_x|$ is understood in the Sobolev sense.

\begin{lemma}[Uniform continuity and limits on the line \cite{gallo}]\label{lem:line-continuity}
Every $v\in X^1(\R)$ has a uniformly continuous representative satisfying
\[
 |v(x)-v(y)|\le |x-y|^{1/2}\|v_x\|_2,\qquad x,y\in\R.
\]
If also $v\in L^p(\R)$ for some $1\le p<\infty$, then $v(x)\to0$ as $x\to\pm\infty$.
In particular, every $w\in H^1(\R)$ tends to $0$ at both ends, and every field $1+w\in1+H^1(\R)$ tends to $1$ there.
\end{lemma}
\begin{proof}
The estimate follows by integrating $v_x$ and applying Cauchy--Schwarz; see also \cite[Proposition~2.2]{gallo}. If decay failed, uniform continuity would give infinitely many disjoint intervals of a fixed positive length on which $|v|$ is bounded below by a positive constant, contradicting $v\in L^p$.
For $w\in H^1(\R)$, the one-dimensional Sobolev embedding gives $w\in L^\infty$, so $w\in X^1\cap L^2$. Applying the preceding assertion with $p=2$ gives $w(x)\to0$, and hence $1+w(x)\to1$, at both ends.
\end{proof}

For $v\in\E_1$, the pointwise bound $\bigl||v|-1\bigr|\le \bigl||v|+1\bigr|\cdot \bigl||v|-1\bigr|= \bigl||v|^2-1\bigr|$ gives $|v|-1\in L^2(\R)$. By the Sobolev chain rule, the weak derivative of $|v|$ satisfies $\bigl|\partial_x|v|\bigr|\le |v_x|$ almost everywhere. Thus $|v|-1\in H^1(\R)\hookrightarrow L^\infty(\R)$, so $v\in X^1(\R)$.
The inclusions
\[
 1+H^1(\R)\subsetneq\E_1\subsetneq X^1(\R)
\]
are strict. 
Indeed, if $w\in H^1(\R)$, then $w\in L^2\cap L^\infty$, so $|w|^2\in L^2$ and $|1+w|^2-1=2\operatorname{Re}w+|w|^2\in L^2$. Together with $\partial_x(1+w)=w_x\in L^2$, this proves $1+H^1\subseteq\E_1$. The field $v(x)=e^{i\log(1+x^2)}$ of B\'ethuel et al.\ \cite{bgss-black} belongs to $\E_1$, since $|v|=1$ and $v_x=2ix(1+x^2)^{-1}v\in L^2$. It does not belong to $1+H^1(\R)$: it has no limit at either end, whereas every $1+H^1$ field tends to $1$ at both ends by Lemma~\ref{lem:line-continuity}. This proves that the first inclusion is strict. Finally, the constant $2$ belongs to $X^1\setminus\E_1$.

We use three elementary facts about these spaces. First, the product rule makes $X^1$ and $H^1$ algebras, with multiplication also continuous from $X^1\times H^1$ to $H^1$. Second, if $g$ and its distributional second derivative $g_{xx}$ belong to $L^2(\R)$, Plancherel and Cauchy--Schwarz give
\begin{equation}\label{eq:line-interpolation}
 \|g_x\|_2^2\le\|g\|_2\|g_{xx}\|_2,\qquad
 \|g\|_{H^2}\le C\bigl(\|g\|_2+\|g_{xx}\|_2\bigr).
\end{equation}
Third, every $v\in\E_1$ satisfies $\|v\|_\infty\ge1$, since otherwise $1-|v|^2$ would be bounded below by a positive constant on $\R$.

The following lemma is a version of the one-dimensional Gagliardo--Nirenberg inequality
\[
 \|h\|_\infty^3\le C\|h\|_4^2\|h_x\|_2,
 \qquad h\in L^4(\R),\quad h_x\in L^2(\R),
\]
applied to the excess amplitude $h=(|v|-A)_+$.

\begin{lemma}[Estimation of excess amplitude]\label{lem:line-amplitude}
Let $v\in X^1(\R)$ and let $b$ be a measurable density with $0\le b\le A^2$ for some $A\ge0$. If $|v|^2-b\in L^2(\R)$, then
\[
 \big\|(|v|-A)_+\big\|_\infty^3\le 3\big\||v|^2-b\big\|_2\,\|v_x\|_2 .
\]
\end{lemma}
\begin{proof}
Put $h=(|v|-A)_+$. Where $|v|>A$ we have $(|v|-A)^2\le |v|^2-A^2\le |v|^2-b$. Thus $\|h\|_4^2\le\||v|^2-b\|_2$, and the chain rule gives $|h_x|\le |v_x|$ almost everywhere. Since $h\in X^1\cap L^4$, Lemma~\ref{lem:line-continuity} gives $h(x)\to0$ as $x\to-\infty$. Integrating the derivative of $h^3$ from $-\infty$ to $x$ gives
\[
 h(x)^3\le 3\int_{\R}h^2|h_x|\dd\le 3\|h\|_4^2\|v_x\|_2 \le3\big\||v|^2-b\big\|_2\,\|v_x\|_2.
\]
\end{proof}

\begin{lemma}[Free increments {\cite{gallo,gerard}}]\label{lem:free}
Let $S(t)=e^{it\partial_x^2}$. For $a\in X^1(\R)$, $(S(t)-I)a\in C(\R;H^1(\R))$, and
\[
 \|(S(t)-I)a\|_2\le C|t|^{1/2}\|a_x\|_2,\qquad
 \partial_x(S(t)-I)a=(S(t)-I)a_x.
\]
In particular $S(t)$ is a strongly continuous group on $X^1(\R)$.
\end{lemma}
\begin{proof}
The group constructed by Gallo \cite[Section~3]{gallo} agrees with the Fourier multiplier $e^{-it\xi^2}$ in tempered distributions.
Define the smooth multiplier $m_t(\xi)=(e^{-it\xi^2}-1)/\xi^2$, with $m_t(0)=-it$. Then $(S(t)-I)a=-m_t(D)\partial_xa_x$ in tempered distributions. The bound $\sup_\xi |\xi|\,|m_t(\xi)|\le C\sqrt{|t|}$ gives the $L^2$ estimate by Plancherel. The derivative identity and strong continuity on $L^2$ give $H^1$ continuity at zero, and the group law gives it at every time. Since $H^1(\R)\hookrightarrow L^\infty(\R)$, this also gives strong continuity in $X^1$. See \cite[Lemma~3]{gerard} and \cite[Section~3]{gallo}.
\end{proof}

\begin{lemma}[An elementary bound for the free group]\label{lem:free-group-bound}
There is a universal $C$ such that, for every $T>0$ and $a\in X^1(\R)$,
\begin{equation}\label{eq:free-group-bound}
 \|S(s)a\|_{X^1}\le K_T\|a\|_{X^1},\qquad
 K_T=1+CT^{1/4},\qquad |s|\le T.
\end{equation}
\end{lemma}
\begin{proof}
Lemma~\ref{lem:free}, the Gagliardo--Nirenberg inequality and the $L^2$ isometry of $S(s)$ give
\[
 \begin{aligned}
 \|(S(s)-I)a\|_\infty
 &\le C\|(S(s)-I)a\|_2^{1/2}\|(S(s)-I)a_x\|_2^{1/2}\\
 &\le C|s|^{1/4}\|a_x\|_2.
 \end{aligned}
\]
Adding $\|a\|_\infty$ and using $\|\partial_xS(s)a\|_2=\|a_x\|_2$ proves \eqref{eq:free-group-bound}.
\end{proof}
A sharper large-time estimate, recovering Gallo's sharp exponent $1/6$ \cite[Theorem~1.1]{gallo-growth}, follows from Remark~\ref{rem:cubic-amplitude}. 

We recall two standard facts about differentiation of norms.

\begin{lemma}[Differentiation of $L^2$ norms]\label{lem:norm-differentiation}
Let $I$ be an interval containing $0$ and $S(t)=e^{it\partial_x^2}$.

\smallskip\noindent\textup{(i)} Let $v_0\in L^2(\R)$, $G\in C(I;L^2(\R))$, and
\[
 v(t)=S(t)v_0+\int_0^t S(t-s)G(s)\,ds.
\]
Then $t\mapsto\|v(t)\|_2^2$ is continuously differentiable, with
\[
 \tfrac12\frac{d}{dt}\|v(t)\|_2^2
 =\operatorname{Re}\int_{\R}\overline{v(t)}\,G(t)\dd .
\]
This identity requires no a priori assumption that $v_t(t)\in L^2(\R)$.

\smallskip\noindent\textup{(ii)} Let $r,b:I\to[0,\infty)$ with $r^2\in C^1(I)$ and $b\in C(I)$.
If
\[
 \frac12\left|\frac{d}{dt}\bigl(r(t)^2\bigr)\right|
 \le r(t)b(t) \qquad \text{for almost every }t\in I,
\]
then
\[
 |r(t)-r(s)|\le\int_s^t b(\tau)\,d\tau,\qquad s,t\in I,\quad s\le t.
\]
In particular, $r$ is locally Lipschitz and $|r'|\le b$ almost everywhere.
\end{lemma}

Finally, note that the equivalence between writing the nonlinearity as $q(|u|^2-1)$ or $q|\psi|^2$ holds in our setting. More generally,

\begin{lemma}[Phase rotation for density responses on the line]\label{lem:phase-rotation}
Let $I$ be an interval containing $0$, let $\omega\in\R$, and let $f$ be any of the density responses covered by Propositions~\ref{prop:bounded-slope}--\ref{prop:subcubic-powers}. Suppose $u\in C(I;\X)$ is a mild solution of the rotated equation
\[
 iu_t+u_{xx}+\bigl(f(|u|^2)-\omega\bigr)u=0 .
\]
Then $\psi(t)=e^{i\omega t}u(t)$ is a mild solution of the unrotated equation \eqref{eq:nls-unrotated} on $I$. The rotation is invertible and preserves the initial datum, maximal lifespan, density, mass current and $L^2$ norm of the spatial derivative.
\end{lemma}
\begin{proof}
In all the nonlinearities covered here, for every $u\in C(I;X^1)$ it follows that $(f(|u|^2)-\omega)u\in C(I;X^1)$. The mild solution in the rotated frame satisfies
\[
 S(-t)u(t)=u(0)+i\int_0^t S(-s)\bigl[(f(|u(s)|^2)-\omega)u(s)\bigr]\,ds .
\]
The integrand is continuous in $X^1$ and thus  $t\mapsto S(-t)u(t)$ belongs to $C^1(I;X^1)$. Set $\psi(t)=e^{i\omega t}u(t)$. Integrating $\frac{d}{dt} (e^{i\omega t} S(-t)u(t))$ from $0$ to $t$ gives
\[
 \begin{aligned}
 S(-t)\psi(t)
 &=u(0)+i\int_0^t S(-s)\bigl[\omega\psi(s)+(f(|u(s)|^2)-\omega)\psi(s)\bigr]\,ds\\
 &=\psi(0)+i\int_0^t S(-s)\bigl[f(|\psi(s)|^2)\psi(s)\bigr]\,ds,
 \end{aligned}
\]
since $|\psi|=|u|$. This is the Duhamel formula for \eqref{eq:nls-unrotated}. The inverse rotation gives the converse on every interval of existence. Finally, $\psi_x=e^{i\omega t}u_x$ and $\bar\psi\psi_x=\bar u u_x$.
\end{proof}

\begin{lemma}[Composition estimates for density responses]\label{lem:density-composition}
Let $f\in C^2([0,\infty);\R)$ and let $n,b\in X^1$ be nonnegative with $\|n\|_{X^1},\|b\|_{X^1}\le R$. Then
\begin{equation}\label{eq:density-composition}
 \|f(n)-f(b)\|_{X^1}\le C_{f,R}\|n-b\|_{X^1}.
\end{equation}
If in addition $n-b\in L^2$, then $n-b\in H^1$ and
\[
 \|f(n)-f(b)\|_{H^1}\le C_{f,R}\|n-b\|_{H^1}.
\]
\end{lemma}

\begin{proof}
Put $M_k=\sup_{[0,R]}|f^{(k)}|$ for $k=1,2$. The mean value theorem gives $|f(n)-f(b)|\le M_1|n-b|$, and the Sobolev chain rule gives
\[
 \partial_x(f(n)-f(b))=f'(n)(n_x-b_x)+(f'(n)-f'(b))b_x.
\]
Thus the derivative has $L^2$ norm at most $M_1\|n_x-b_x\|_2+M_2R\|n-b\|_\infty$. This proves the $X^1$ estimate. For the $H^1$ estimate, take the $L^2$ norm of the pointwise bound and use $H^1(\R)\hookrightarrow L^\infty(\R)$ in the derivative bound.
\end{proof}

%\begin{remark}[The cubic NLS]
%With $f(s)=qs$ and $\omega=q$, Lemma~\ref{lem:phase-rotation} identifies the global mild solution $u$ of \eqref{eq:gp} with the solution $\psi(t)=e^{iqt}u(t)$ of the focusing cubic nonlinear Schr\"odinger equation
%\begin{equation}\label{eq:cubic}
% i\psi_t+\psi_{xx}+q|\psi|^2\psi=0,\qquad \psi(0)=u_0,\qquad q>0.
%\end{equation}
%For $u_0\in 1+H^1(\R)$, the corresponding perturbation is $\psi(t)-e^{iqt}=e^{iqt}(u(t)-1)\in H^1(\R)$.
%
%\end{remark}

\section{Proof of Theorem~\ref{thm:main}}\label{sec:cubic}
We first extract some technical justifications that can be reused in the proofs of Propositions~\ref{prop:bounded-slope}--\ref{prop:subcubic-powers}. 

\begin{lemma}[Density identity; persistence and Lipschitz estimates in the energy and affine spaces]\label{lem:restricted-flows}
Consider the equation
\begin{equation}\label{eq:restricted-flow}
 iu_t+u_{xx}+\mathcal N(u)=0,\qquad u(0)=u_0.
\end{equation}
Assume that $\mathcal N:\mathbb C\to\mathbb C$ has the following properties.
\begin{itemize}
\item For every $B>0$ there is a constant $L_B$ such that
      $|\mathcal N(z)-\mathcal N(w)|\le L_B|z-w|$ whenever $|z|,|w|\le B$.
\item $\mathcal N(z)=0$ whenever $|z|=1$.
\item $\operatorname{Im}(\bar z\,\mathcal N(z))=0$ for every $z\in\mathbb C$.
\end{itemize}
Let $J$ be a compact interval containing $0$.

\smallskip\noindent\textup{(i) Density identity.} Every mild solution $u\in C(J;X^1)$ of \eqref{eq:restricted-flow} satisfies, for each $t\in J$,
\begin{equation}\label{eq:restricted-density-identity}
 |u(t)|^2-|u_0|^2=-2\partial_x\int_0^t\operatorname{Im}(\bar u u_x)(s)\,ds.
\end{equation}

For parts~\textup{(ii)} and~\textup{(iii)}, assume moreover that \eqref{eq:restricted-flow} is locally well-posed in $X^1$ and that the mild solution map
\[
 \Phi\colon X^1\longrightarrow C(J;X^1),\qquad u_0\longmapsto u(\cdot),
\]
is defined, bounded and Lipschitz on each initial $X^1$ ball.

\smallskip\noindent\textup{(ii) Energy space.} If $u_0\in\E_1$, then $u\in C(J;\E_1)$ and the solution map is Lipschitz on bounded sets in the energy metric $d_{\E_1}$, uniformly on $J$.

\smallskip\noindent\textup{(iii) Affine space.} If $u_0\in1+H^1$, then $u-1\in C(J;H^1)$ and the solution map is Lipschitz on bounded sets in the affine metric $\|u_1-u_2\|_{H^1}$, uniformly on $J$.
\end{lemma}
\begin{proof}
Note that $v\in X^1$ implies   $\mathcal N(v)\in X^1;$ thus
\[
 S(-t)u(t)=u_0+i\int_0^t S(-s)\mathcal N(u(s))\,ds 
\]
makes sense with $S(t)\colon X^1\to X^1$ as in Lemmas~\ref{lem:free} and~\ref{lem:free-group-bound}.

\smallskip\noindent\textup{(i)} We first derive the density identity for a single mild solution. Let $u_\varepsilon$ and $N_\varepsilon$ be the spatial mollifications of $u$ and $\mathcal N(u)$, respectively. Mollifying the equation and taking the real part after multiplication by $2i\bar u_\varepsilon$ gives
\[
 \partial_t|u_\varepsilon|^2
 +2\partial_x\operatorname{Im}(\bar u_\varepsilon\partial_xu_\varepsilon)
 =-2\operatorname{Im}(\bar u_\varepsilon N_\varepsilon).
\]
Here $u_\varepsilon\to u$ and $N_\varepsilon\to\mathcal N(u)$ locally uniformly, and $\partial_xu_\varepsilon\to u_x$ in $C(J;L^2)$. Thus the right-hand side tends locally to $-2\operatorname{Im}(\bar u\mathcal N(u))=0$. Passing to the limit in distributions gives $\partial_t|u|^2=-2\partial_x\operatorname{Im}(\bar u u_x)$. Testing against a spatial test function and integrating from $0$ to $t$ gives \eqref{eq:restricted-density-identity}.

\smallskip\noindent\textup{(ii)} Now take two such solutions $u_1,u_2$, choose $B\ge1$ with $\sup_{s\in J}\|u_i(s)\|_{X^1}\le B$ for $i=1,2$, and denote $\rho=|u_1|^2-|u_2|^2$. Set $A=\bar u_1\partial_xu_1-\bar u_2\partial_xu_2$. Adding and subtracting $\bar u_2\partial_xu_1$ gives
\[
 \|A(s)\|_2\le \|\partial_xu_1(s)\|_2\|u_1(s)-u_2(s)\|_\infty
 +\|u_2(s)\|_\infty\|\partial_x(u_1-u_2)(s)\|_2
 \le B\|u_1(s)-u_2(s)\|_{X^1}.
\]
For $a,b\in J$ with $a\le b$, denote the difference of the net mass transfers over $[a,b]$ by
\[
 G:=-2\int_a^b\operatorname{Im} A(s)\,ds
 =-2\int_a^b\bigl[\operatorname{Im}(\bar u_1\partial_xu_1)
                  -\operatorname{Im}(\bar u_2\partial_xu_2)\bigr](s)\,ds .
\]
Applying the density identity \eqref{eq:restricted-density-identity} to each solution gives $G_x=\rho(b)-\rho(a)$. Since $\rho_x=2\operatorname{Re} A$, it follows that $G_{xx}=2\operatorname{Re}(A(b)-A(a))$. Thus
\[
 \|G\|_2\le 2(b-a)\sup_{[a,b]}\|A\|_2,\qquad
 \|G_{xx}\|_2\le 4\sup_{[a,b]}\|A\|_2.
\]
In particular, $G\in H^2$ by \eqref{eq:line-interpolation}, so $G$ vanishes at $-\infty$ and $G(x)=\lim_{y\to-\infty}\int_y^x[\rho(b,\xi)-\rho(a,\xi)]\,d\xi$.
Substituting these bounds in $\|G_x\|_2^2\le\|G\|_2\|G_{xx}\|_2$ from \eqref{eq:line-interpolation} gives
\begin{equation}\label{eq:energy-density-difference}
 \|\rho(b)-\rho(a)\|_2
 \le 2\sqrt{2}\,B\sqrt{b-a}\sup_{a\le s\le b}\|u_1(s)-u_2(s)\|_{X^1}.
\end{equation} 

Now take $u_2\equiv1$, which satisfies \eqref{eq:restricted-density-identity} since both sides vanish. Then $\rho^1=|u_1|^2-1$. For $u_1(0)\in\E_1$ we have $\rho^1(0)\in L^2$, so \eqref{eq:energy-density-difference} gives $\rho^1\in C(J;L^2)$. Together with $u_1\in C(J;X^1)$, this gives $u_1\in C(J;\E_1)$ by \eqref{eq:energy-space}.

For Lipschitz dependence, take two initial data $u_1(0),u_2(0)$ in a fixed $\E_1$ ball and let $u_1,u_2$ be their solutions. The initial ball is bounded in $X^1$ because $\|v-w\|_{X^1}\le d_{\E_1}(v,w)$. Thus the assumed bound on the $X^1$ flow gives a common $B$ for these solutions on $J$, and its Lipschitz estimate gives a constant $K\ge1$ such that
\[
 \sup_{s\in J}\|u_1(s)-u_2(s)\|_{X^1}
 \le K\|u_1(0)-u_2(0)\|_{X^1}.
\]
It remains to bound the density difference $\rho(t)=|u_1(t)|^2-|u_2(t)|^2=(|u_1(t)|^2-1)-(|u_2(t)|^2-1)$. Its initial value $\rho(0)$ belongs to $L^2$ since it is the difference of the two initial density defects. Applying \eqref{eq:energy-density-difference} between $0$ and $t$ in either time direction gives
\[
 \|\rho(t)\|_2
\le \|\rho(0)\|_2+2\sqrt{2}\,B\sqrt{|t|}\,K\|u_1(0)-u_2(0)\|_{X^1}.
\]
Adding the $X^1$ estimate to this density estimate recovers all three terms of the energy metric \eqref{eq:energy-space}. With $T=\max_{s\in J}|s|$, we obtain
\[
 \begin{aligned}
 d_{\E_1}(u_1(t),u_2(t))
 &=\|u_1(t)-u_2(t)\|_{X^1}+\|\rho(t)\|_2\\
 &\le K(1+2\sqrt{2}\,B\sqrt T)\|u_1(0)-u_2(0)\|_{X^1}+\|\rho(0)\|_2\\
 &\le K(1+2\sqrt{2}\,B\sqrt T)\,d_{\E_1}(u_1(0),u_2(0)).
 \end{aligned}
\]
%The last inequality uses $K\ge1$ and the same metric formula at $t=0$. This proves Lipschitz dependence on the initial energy ball, uniformly for $t\in J$.

\smallskip\noindent\textup{(iii)} Suppose now that $u_i(0)-1\in H^1$. Since $1+H^1\subset\E_1$, part~\textup{(ii)} gives $|u_i|^2-1\in C(J;L^2)$. To prove $u_i-1\in C(J;H^1)$, it remains to show $u_i-1\in C(J;L^2)$, since $\partial_xu_i\in C(J;L^2)$ already follows from $u_i\in C(J;X^1)$. We use the density bound to control the forcing term in the Duhamel formula for $u_i-1$.

Take a bound $B$ of the $X^1$-ball where our solution lives, as above. 
Let $L_B$ be a Lipschitz constant for $\mathcal N$ on $|z|\le B$. For $0<|z|\le B$, compare $\mathcal N(z)$ with $\mathcal N(z/|z|)=0$. The Lipschitz bound gives
\[
 |\mathcal N(z)|\le L_B\left|z-\frac{z}{|z|}\right|
 = L_B\frac{\big||z|^2-1\big|}{|z|+1}
 \le L_B\big||z|^2-1\big|.
\]
At $z=0$ the same bound follows by comparison with $\mathcal N(1)=0$. Thus 
$$
\|\mathcal N(u_i(t))\|_2\le L_B\||u_i(t)|^2-1\|_2
$$
is bounded on $J$.  Since $S(t)1=1$, Duhamel's formula gives
\[
 u_i(t)-1=S(t)(u_i(0)-1)+i\int_0^t S(t-s)\mathcal N(u_i(s))\,ds.
\]
The right-hand side is continuous in $L^2$. Together with $\partial_xu_i\in C(J;L^2)$ this proves affine persistence.

For Lipschitz dependence, write $w=u_1-u_2$ and take the difference of the two Duhamel formulas. Unitarity of $S(t)$ on $L^2$ and the pointwise bound $|\mathcal N(u_1)-\mathcal N(u_2)|\le L_B|w|$ give, for $t\ge0$,
\[
 \|w(t)\|_2\le\|w(0)\|_2+L_B\int_0^t\|w(s)\|_2\,ds .
\]
Gronwall's inequality, also applied backward in time, gives $\|w(t)\|_2\le e^{L_B|t|}\|w(0)\|_2$. To control the derivative, fix an affine initial ball. Since $H^1\hookrightarrow X^1$, the initial data in this ball are bounded in $X^1$. The $X^1$ flow therefore supplies common constants $B$ and $K$ as in~\textup{(ii)}, and
\[
 \|w_x(t)\|_2\le\|w(t)\|_{X^1}\le K\|w(0)\|_{X^1}
 \le K C_{\rm emb}\|w(0)\|_{H^1},
\]
where $C_{\rm emb}$ is the embedding constant. Adding the squares of the $L^2$ and derivative bounds gives
\[
 \sup_{t\in J}\|w(t)\|_{H^1}^2
 \le \bigl(e^{2L_BT}+K^2C_{\rm emb}^2\bigr)\|w(0)\|_{H^1}^2,
 \qquad T=\max_{s\in J}|s|.
\]
\end{proof}

Now, for the proof of Theorem \ref{thm:main}, we start from a local mild solution in $X^1.$ Indeed,
fix $q>0$ and $u_0\in X^1(\R)$; then Gallo's local theory gives a unique maximal mild solution $u\in C(I;X^1)$ of \eqref{eq:gp} on $I=(T_-,T_+)\ni0$ \cite[Theorems~4.1--4.2 and Proposition~4.1]{gallo}. For data in $\{\|u_0\|_{X^1}\le R\}$ there is a common existence time $\tau_R>0$. Since the cubic is Lipschitz on $X^1$ balls, the Duhamel estimate gives Lipschitz dependence in $C([-\tau_R,\tau_R];X^1)$ after shrinking $\tau_R$ if necessary. At a finite endpoint of $I$ the blow-up alternative requires
\[
 \|u(t)\|_{X^1}=\|u(t)\|_\infty+\|u_x(t)\|_2\longrightarrow\infty .
\]
We will show that $\|u(t)\|_{X^1}$ cannot blow up in finite time, thus $I=\mathbb{R}$.

\begin{proof}[\normalfont\bfseries Proof of Theorem~\ref{thm:main}]
The proof is organized in a series of distinct steps. Step~1 constructs the net mass transfer function and establishes the mass-weighted cancellation \eqref{eq:mass-weighted-cancellation} for the local mild solution $u\in C(I;X^1).$  Step~2 establishes the identity for the modified derivative, and Step~3 closes the estimate for $F$ and recovers the kinetic bound. Step~4 controls the amplitude and extends the $X^1$ flow globally. Step~5 proves Lipschitz dependence in $X^1$. Step~6 treats the energy and affine spaces.

\textbf{Step 1. Net mass transfer and the mass-weighted cancellation.}
All limiting arguments below are carried out on compact subintervals of $I$ containing zero. Write $n=|u|^2$ and $n_0=|u_0|^2$. The products
\[
 j=2\operatorname{Im}(\bar u u_x),\qquad n_x=2\operatorname{Re}(\bar u u_x)
 \quad\text{belong to }C(I;L^2).
\]
Gallo's mild solution belongs to $C(J;X^1)$ on every compact interval $J\subset I$ containing $0$. The cubic $\mathcal N(z)=q(|z|^2-1)z$ is Lipschitz on bounded disks, vanishes on the unit circle and satisfies $\bar z\mathcal N(z)\in\R$. Thus Lemma~\ref{lem:restricted-flows}(i) gives \eqref{eq:restricted-density-identity}. Defining $F$ as the time integral of $-j$ in $L^2$, we obtain
\begin{equation}\label{eq:cubic-flux}
 F(t)=-\int_0^t j(s)\,ds,\qquad F_t=-j,\qquad F_x(t)=n(t)-n_0 .
\end{equation}
Since $F\in C^1(I;L^2)$ and $F_{xx}=n_x-\partial_xn_0\in C(I;L^2)$, \eqref{eq:line-interpolation} gives $F\in C(I;H^2)$ and $F_x\in C(I;H^1)$. The squared density increment is therefore
\begin{equation}\label{eq:squared-density-increment}
 \|n(t)-n_0\|_2^2=\|F_x(t)\|_2^2.
\end{equation}

By Lemma~\ref{lem:line-continuity}, $F(t,\cdot)$ has a continuous representative vanishing at both spatial infinities, and
\[
 \int_a^b(n(t,x)-n_0(x))\dd=F(t,b)-F(t,a),\qquad
 F(t,x)=\lim_{a\to-\infty}\int_a^x(n(t,y)-n_0(y))\,dy.
\]
The improper integral need not converge absolutely.
We note that B\'ethuel, Gravejat, Saut and Smets \cite{bgss-kdv2} use a similar device to keep track of net mass transfer in a different context.

Since $F^3\in W^{1,1}(\R)$, with $\|F_xF^2\|_1\le\|F\|_\infty\|F\|_2\|F_x\|_2$, its derivative integrates to zero. This proves the mass-weighted cancellation \eqref{eq:mass-weighted-cancellation}, together with the bound
\begin{equation}\label{eq:cubic-weight}
 \int_{\R}(n-n_0)F^2\dd =\tfrac13\int_{\R}(F^3)_x\dd=0 \quad\Longrightarrow\quad
 \|F(t)u(t)\|_2=\|F(t)u_0\|_2\le \|u_0\|_\infty \|F(t)\|_2.
\end{equation}

\medskip\noindent\textbf{Step 2. The modified derivative and its exact identity.}
We first record the differentiation rule for the squared density increment.

\begin{claim}\label{cl:gradient-norm}
Let $I$ be an open interval and let $F\in C(I;H^2(\R))\cap C^1(I;L^2(\R))$ be real-valued. Then $t\mapsto\|F_x(t)\|_2^2$ is continuously differentiable and
\[
 \tfrac12\frac{d}{dt}\|F_x(t)\|_2^2=-\int_{\R} F_{xx}(t)F_t(t)\dd .
\]
\end{claim}
\noindent\emph{Proof of Claim~\ref{cl:gradient-norm}.}
Fix $t\in I$ and $h\ne0$ with $t+h\in I$. Then
\[
 \begin{aligned}
 \frac{\|F_x(t+h)\|_2^2-\|F_x(t)\|_2^2}{2h}
 &=\int_{\R}\frac{F_x(t+h)+F_x(t)}{2}\,
             \frac{F_x(t+h)-F_x(t)}{h}\dd\\
 &=-\int_{\R}\frac{F_{xx}(t+h)+F_{xx}(t)}{2}\,
             \frac{F(t+h)-F(t)}{h}\dd .
 \end{aligned}
\]
As $h\to0$, the two factors converge in $L^2$ to $F_{xx}(t)$ and $F_t(t)$, respectively. Cauchy--Schwarz gives the  limit, and its continuity follows from $F_{xx},F_t\in C(I;L^2)$. The proof of Claim~\ref{cl:gradient-norm} is complete.

Define
\[
 Y=u_x+\tfrac q2Fu\in C(I;L^2).
\]
Since $F$ is real, $j=2\operatorname{Im}(\bar uY)$. 

\begin{claim*}
The modified derivative satisfies
\begin{equation}\label{eq:cubic-square}
 \|Y(t)\|_2
 =\big\|\partial_xu_0-\tfrac q2F(t)u_0\big\|_2 .
\end{equation}
\end{claim*}

\noindent\emph{Proof of the claim.}
With $S(t)=e^{it\partial_x^2}$, the spatial derivative satisfies
\[
 u_x(t)=S(t)\partial_xu_0
 +iq\int_0^t S(t-s)\partial_x\bigl[(n(s)-1)u(s)\bigr]\,ds .
\]
Since $\partial_x((n-1)u)=n_xu+(n-1)u_x\in C(I;L^2)$, Lemma~\ref{lem:norm-differentiation}(i) applied to $v=u_x$ gives
\[
 \frac{d}{dt}\|u_x(t)\|_2^2
 =2\operatorname{Re}\int iq\bar u_x\bigl(n_xu+(n-1)u_x\bigr)\dd
 =q\int n_xj\dd .
\]
Here $iq(n-1)|u_x|^2$ is purely imaginary.
Step~1 gives $F\in C(I;H^2)\cap C^1(I;L^2)$. Substituting $n_x=F_{xx}+\partial_x n_0$ and $j=-F_t$ and applying Claim~\ref{cl:gradient-norm} gives
\[
 \begin{aligned}
 \frac{d}{dt}\|u_x(t)\|_2^2
 &=-q\int(F_{xx}+\partial_x n_0)F_t\dd\\
 &=\tfrac q2\frac{d}{dt}\|F_x(t)\|_2^2
 -q\frac{d}{dt}\int F\,\partial_x n_0\dd .
 \end{aligned}
\]
Integrating from zero to $t$, using $F(0)=0$, and then integrating by parts with $F_{xx}=n_x-\partial_xn_0$, we obtain the kinetic energy increment
\[
 \begin{aligned}
 \|u_x(t)\|_2^2-\|\partial_xu_0\|_2^2
 &=\tfrac q2\|F_x\|_2^2-q\int F\,\partial_x n_0\dd\\
 &=\tfrac q2\int F_x(n(t)-n_0)\dd-q\int F\,\partial_x n_0\dd\\
 &=-\tfrac q2\int F(n(t)-n_0)_x\dd-q\int F\,\partial_x n_0\dd\\
 &=-\tfrac q2\int F(n_x+\partial_x n_0)\dd .
 \end{aligned}
\]
Here $F_x=n(t)-n_0$, and the integration by parts is justified by $F\in H^2(\R)$, as established in Step~1. (In particular, we manage to bound the increment of the kinetic energy in terms of $F$ without ever requiring the potential energy to be finite.)

So finally
\begin{equation}\label{eq:balance}
	 \|u_x(t)\|_2^2 + \tfrac q2\int F n_x \dd = \|\partial_xu_0\|_2^2 -\tfrac q2\int F \partial_x n_0\dd.
\end{equation}
To conclude, expanding the square and using \eqref{eq:balance} and the mass-weighted cancellation \eqref{eq:cubic-weight} gives
\[
 \begin{aligned}
 \|Y(t)\|_2^2
 &=\|u_x(t)\|_2^2+\tfrac q2\int Fn_x\dd+\tfrac{q^2}{4}\|Fu\|_2^2\\
 &=\|\partial_xu_0\|_2^2-\tfrac q2\int F\,\partial_xn_0\dd
   +\tfrac{q^2}{4}\|Fu_0\|_2^2\\
 &=\|\partial_xu_0\|_2^2-q\operatorname{Re}\int F\overline{u_0}\,\partial_xu_0\dd
   +\tfrac{q^2}{4}\|Fu_0\|_2^2\\
 &=\big\|\partial_xu_0-\tfrac q2Fu_0\big\|_2^2 .
 \end{aligned}
\]
This proves the claim.

In particular,
\begin{equation}\label{eq:cubic-Y-bound}
 \|Y(t)\|_2\le \|\partial_xu_0\|_2+\tfrac q2\|u_0\|_\infty \|F(t)\|_2.
\end{equation}

\medskip\noindent\textbf{Step 3. Close the estimate for $F$ and recover the kinetic bound.}
%Since $F\in C^1(I;L^2)$, its norm $\|F(t)\|_2$ is locally Lipschitz. 
Using $F_t=-j$, 
\[
 \begin{aligned}
 \tfrac12\left|\frac{d}{dt}\|F(t)\|_2^2\right|
 &=\left|\int Fj\dd\right|
 \le 2\|Fu\|_2\|Y\|_2\\
 &\le 2\|Fu_0\|_2\bigl(\|\partial_xu_0\|_2+\tfrac q2\|u_0\|_\infty \|F(t)\|_2\bigr),
 \end{aligned}
\]
where we used $j=2\operatorname{Im}(\bar uY)$, the mass-weighted cancellation \eqref{eq:cubic-weight}, and the bound \eqref{eq:cubic-Y-bound} for $Y$.

Using $\|Fu_0\|_2\le\|u_0\|_\infty\|F(t)\|_2$, we apply Lemma~\ref{lem:norm-differentiation}(ii) with $r(t)=\|F(t)\|_2$ to obtain
\[
 |r'(t)|\le 2\|u_0\|_\infty\|\partial_xu_0\|_2+q\|u_0\|_\infty^2r(t)
 \quad\text{a.e.}
\]
Since $r(0)=0$, Gronwall's inequality in either time direction gives
\[
 \|F(t)\|_2 \le 2\|u_0\|_\infty \|\partial_xu_0\|_2\int_0^{|t|}e^{q\|u_0\|_\infty^2s}\,ds.
\]
For the kinetic norm, the definition of $Y$ and \eqref{eq:cubic-square}, together with the reverse triangle inequality, give
\[
 \begin{aligned}
 Y=u_x+\tfrac q2Fu
 &\quad\Longrightarrow\quad
 \bigl|\|u_x\|_2-\|Y\|_2\bigr|\le\tfrac q2\|Fu\|_2,\\
 \|Y\|_2=\bigl\|\partial_xu_0-\tfrac q2Fu_0\bigr\|_2
 &\quad\Longrightarrow\quad
 \bigl|\|Y\|_2-\|\partial_xu_0\|_2\bigr|\le\tfrac q2\|Fu_0\|_2.
 \end{aligned}
\]
Adding these inequalities and using the mass-weighted cancellation \eqref{eq:cubic-weight} yields
\begin{equation}\label{eq:cubic-gradient}
 \bigl|\|u_x(t)\|_2-\|\partial_xu_0\|_2\bigr|\le  q\|u_0\|_\infty \|F(t)\|_2.
\end{equation}
Thus, 
\begin{equation}\label{eq:cubic-current-growth}
 \begin{aligned}
 \|u_x(t)\|_2&\le \|\partial_xu_0\|_2+q\|u_0\|_\infty \|F(t)\|_2\le \|\partial_xu_0\|_2(2e^{q\|u_0\|_\infty^2|t|}-1).
 \end{aligned}
\end{equation}
This proves the kinetic norm bound \eqref{eq:cubic-kinetic-growth} throughout the local lifespan.

\medskip\noindent\textbf{Step 4. Recover the amplitude and continue the solution.}
Set $\delta(t)=\|n(t)-n_0\|_2=\|F_x(t)\|_2$. The integration by parts in Step~2 and the mass-weighted cancellation \eqref{eq:cubic-weight} give
\begin{equation}\label{eq:cubic-density-growth}
 \begin{aligned}
 \delta(t)^2=\|F_x\|_2^2
 &=-\int F(n_x-\partial_x n_0)\dd\\
 &\le 2\bigl(\|Fu\|_2\,\|u_x(t)\|_2+\|Fu_0\|_2\,\|\partial_xu_0\|_2\bigr)\\
 &\le 2\|u_0\|_\infty \|F(t)\|_2\bigl(\|u_x(t)\|_2+\|\partial_xu_0\|_2\bigr).
 \end{aligned}
\end{equation}
Denoting $M=e^{q\|u_0\|_\infty^2|t|}$ for brevity, the bounds in \eqref{eq:cubic-current-growth} give
\[
 \begin{aligned}
 q\|u_0\|_\infty\|F(t)\|_2
 &\le 2\|\partial_xu_0\|_2(M-1),\\
 \|u_x(t)\|_2+\|\partial_xu_0\|_2
 &\le \|\partial_xu_0\|_2(2M-1)+\|\partial_xu_0\|_2
 =2\|\partial_xu_0\|_2M.
 \end{aligned}
\]
Plugging this into \eqref{eq:cubic-density-growth} yields
\begin{equation}\label{eq:cubic-density-explicit}
 \begin{aligned}
 q\delta(t)^2
 &\le 2\bigl(q\|u_0\|_\infty\|F(t)\|_2\bigr)
       \bigl(\|u_x(t)\|_2+\|\partial_xu_0\|_2\bigr)\\
 &\le 2\bigl(2\|\partial_xu_0\|_2(M-1)\bigr)
       \bigl(2\|\partial_xu_0\|_2M\bigr)\\
 &=8\|\partial_xu_0\|_2^2M(M-1).
 \end{aligned}
\end{equation}
At any point where $|u(t,x)|>\|u_0\|_\infty$, we have
\[
 \begin{aligned}
 \bigl(|u(t,x)|-\|u_0\|_\infty\bigr)^2
 &\le \bigl(|u(t,x)|-\|u_0\|_\infty\bigr)
        \bigl(|u(t,x)|+\|u_0\|_\infty\bigr)\\
 &=|u(t,x)|^2-\|u_0\|_\infty^2\\
 &\le |u(t,x)|^2-|u_0(x)|^2\\
 &=n(t,x)-n_0(x).
 \end{aligned}
\]
At all other points the excess amplitude is zero. Taking $L^2$ norms therefore yields
\[
 \big\|(|u(t)|-\|u_0\|_\infty)_+\big\|_4^2
 \le \big\|(n(t)-n_0)_+\big\|_2
 \le \|n(t)-n_0\|_2=\delta(t).
\]
Lemma~\ref{lem:line-amplitude} now applies to give
\[
 \big\|(|u(t)|-\|u_0\|_\infty)_+\big\|_\infty^3\le 3\delta(t)\|u_x(t)\|_2
 \le 12\sqrt{2}\,q^{-1/2}\|\partial_xu_0\|_2^2 M^2 .
\]
Taking cube roots and adding $\|u_0\|_\infty$ leads to the amplitude bound \eqref{eq:cubic-amplitude-growth}. 

Both terms of the $X^1$ norm are now bounded on $I\cap[-T,T]$, uniformly for $u_0$ in any fixed $X^1$ ball. The blow-up alternative gives $I=\R$. For $d_{X^1}(u_0,1)<R$ we have $\|u_0\|_\infty\le 1+R$ and $\|\partial_xu_0\|_2\le R$, so the bounds \eqref{eq:cubic-kinetic-growth} and \eqref{eq:cubic-amplitude-growth} give \eqref{eq:cubic-metric-growth} in $X^1$.

\begin{remark}[Amplitude growth and the free Schr\"odinger group]\label{rem:cubic-amplitude}

The amplitude bound can be kept finite at $q=0$. Recall the notation
\[
 \begin{aligned}
 \delta(t)&=\big\||u(t)|^2-|u_0|^2\big\|_2,\qquad
 M=e^{q\|u_0\|_\infty^2|t|},\\
 F(t)&=-2\int_0^t\operatorname{Im}\bigl(\overline{u(s)}u_x(s)\bigr)\,ds.
 \end{aligned}
\]
The estimates of Steps~3--4 prior to division by $q$ hold also at $q=0$. Using $\int_0^{|t|}e^{q\|u_0\|_\infty^2s}\,ds\le |t|M$, they give
\[
 \begin{aligned}
 \|F(t)\|_2&\le 2\|u_0\|_\infty\|\partial_xu_0\|_2|t|M,\\
 \|u_x(t)\|_2+\|\partial_xu_0\|_2&\le 2\|\partial_xu_0\|_2M,\\
 \delta(t)^2&\le 2\|u_0\|_\infty\|F(t)\|_2
       \bigl(\|u_x(t)\|_2+\|\partial_xu_0\|_2\bigr)\\
 &\le 8\|u_0\|_\infty^2\|\partial_xu_0\|_2^2|t|M^2.
 \end{aligned}
\]
Combining these bounds with Lemma~\ref{lem:line-amplitude} gives
\begin{equation}\label{eq:cubic-amplitude}
 \begin{aligned}
 \big\|(|u(t)|-\|u_0\|_\infty)_+\big\|_\infty^3
 &\le 3\delta(t)\|u_x(t)\|_2\\
 &\le 12\sqrt{2}\,\|u_0\|_\infty \|\partial_xu_0\|_2^2 |t|^{1/2}M^2 .
 \end{aligned}
\end{equation}

Taking the cube root in \eqref{eq:cubic-amplitude} gives, for every $q\ge0$ and $t\in\R$,
\begin{equation}\label{eq:cubic-free-limit}
 \|u(t)\|_\infty\le \|u_0\|_\infty
 + C\|u_0\|_\infty^{1/3}\|\partial_x u_0\|_2^{2/3}|t|^{1/6}
 e^{\frac23 q\|u_0\|_\infty^2|t|},
\end{equation}
with $C=c_*=(12\sqrt{2})^{1/3}$. For $a\in X^1(\R)$, the free solution $S(s)a=e^{is\Delta}a$ satisfies $\|\partial_xS(s)a\|_2=\|a_x\|_2$. Weighted AM--GM gives
\[
 \|a\|_\infty^{1/3}\|a_x\|_2^{2/3}
 \le \tfrac13\|a\|_\infty+\tfrac23\|a_x\|_2
 \le \|a\|_{X^1}.
\]
Combining these facts with \eqref{eq:cubic-free-limit} at $q=0$, we obtain
\[
 \begin{aligned}
 \|S(s)a\|_{X^1}
 &\le \|a\|_\infty+\|a_x\|_2
      +c_*|s|^{1/6}\|a\|_\infty^{1/3}\|a_x\|_2^{2/3}\\
 &\le (1+c_*|s|^{1/6})\|a\|_{X^1}.
 \end{aligned}
\]
This recovers the growth bound with Gallo's sharp exponent $1/6$ \cite[Theorem~1.1]{gallo-growth}.
\end{remark}

\medskip\noindent\textbf{Step 5. Lipschitz dependence in $X^1$.}
Fix $R,T>0$ and let $u(t)=\Phi_tu_0$ and $v(t)=\Phi_tv_0$, with $d_{X^1}(u_0,1),d_{X^1}(v_0,1)<R$. The bounds \eqref{eq:cubic-kinetic-growth} and \eqref{eq:cubic-amplitude-growth} from Step~4 give
\[
\sup_{|t|\le T}\max\{\|u(t)\|_{X^1},\|v(t)\|_{X^1}\}\le B_T,\qquad 
 B_T=\bigl(1+2R+c_*q^{-1/6}R^{2/3}\bigr)e^{q(1+R)^2T}.
\]
Thus a single ball contains both trajectories throughout $[-T,T]$. 

The difference of the two Duhamel formulas for $u(t)$ and $v(t)$ is
\[
 \begin{aligned}
 u(t)-v(t)
 &=S(t)(u_0-v_0)\\
 &\quad+iq\int_0^t S(t-s)
 \bigl[(|u(s)|^2-1)u(s)-(|v(s)|^2-1)v(s)\bigr]\,ds.
 \end{aligned}
\]
Taking $X^1$ norms and applying \eqref{eq:free-group-bound} yields, for $0\le t\le T$,
\[
 \begin{aligned}
 \|u(t)-v(t)\|_{X^1}
 &\le K_T\|u_0-v_0\|_{X^1}\\
 &\quad+qK_T\int_0^t
 \big\|(|u(s)|^2-1)u(s)-(|v(s)|^2-1)v(s)\big\|_{X^1}\,ds\\
 &\le K_T\|u_0-v_0\|_{X^1}
 +qK_T(1+3B_T^2)\int_0^t\|u(s)-v(s)\|_{X^1}\,ds,
 \end{aligned}
\]
where in the last step we used the product rule and the common bound $B_T$:
\[
 \|ab\|_{X^1}\le\|a\|_{X^1}\|b\|_{X^1} \implies
 \big\|(|u|^2-1)u-(|v|^2-1)v\big\|_{X^1}
 \le (1+3B_T^2)\|u-v\|_{X^1}.
\]
Thus, Gronwall's inequality yields
\[
 \sup_{|t|\le T}\|u(t)-v(t)\|_{X^1}
 \le K_Te^{qK_T(1+3B_T^2)T}\|u_0-v_0\|_{X^1}.
\]
This proves \eqref{eq:cubic-metric-dependence} in $X^1$ with an explicit constant depending only on $q,R,T$.

\medskip\noindent\textbf{Step 6. The energy and affine spaces.}
Step~1 verified the hypotheses on the cubic nonlinearity. Steps~4--5 give boundedness and Lipschitz dependence of the $X^1$ flow on initial balls over $[-T,T]$. Parts~\textup{(ii)} and~\textup{(iii)} of Lemma~\ref{lem:restricted-flows} now give persistence and \eqref{eq:cubic-metric-dependence} in both spaces. Uniqueness follows from uniqueness in $X^1$.

For $u_0\in\E_1$, the triangle inequality and \eqref{eq:cubic-density-explicit} give
\[
 \|n(t)-1\|_2\le\|n_0-1\|_2+\delta(t)
 \le \bigl(\|n_0-1\|_2+\sqrt{8/q}\,\|\partial_xu_0\|_2\bigr)M,
 \qquad M=e^{q\|u_0\|_\infty^2|t|}.
\]
This proves \eqref{eq:cubic-energy-density}. Since $d_{\E_1}(u(t),1)=d_{X^1}(u(t),1)+\|n(t)-1\|_2$, adding the $X^1$ growth bound from Step~4 gives \eqref{eq:cubic-metric-growth} on energy balls. We can recover the energy itself. Since $n-1,n_0-1\in H^1$, integration by parts in \eqref{eq:balance} gives
\begin{align*}
 \|u_x(t)\|_2^2-\|\partial_xu_0\|_2^2
 &=\tfrac q2\int F_x(n+n_0-2)\dd\\
 &=\tfrac q2\int (n-n_0)(n+n_0-2)\dd\\
 &=\tfrac q2\bigl(\|n-1\|_2^2-\|n_0-1\|_2^2\bigr).
\end{align*}
Thus $E_q(u(t))=E_q(u_0)$ for all $t$, where $E_q(v)=\|v_x\|_2^2-(q/2)\||v|^2-1\|_2^2$.

For $u_0=1+w_0\in1+H^1$, Lemma~\ref{lem:restricted-flows}(iii) gives $w=u-1\in C(\R;H^1)$. It remains to prove the explicit growth bound \eqref{eq:cubic-affine-growth} for $\|w(t)\|_{H^1}$. Since $(n-1)u=(2\operatorname{Re}w+|w|^2)(1+w)$ is continuous in $L^2$, Lemma~\ref{lem:norm-differentiation}(i) applied to the Duhamel formula for $w$ gives
\[
 \tfrac12\frac{d}{dt}\|w(t)\|_2^2
 =q\operatorname{Re}\int i(n-1)(\bar w+|w|^2)\dd
 =q\operatorname{Re}\int i(n-1)\bar w\dd .
\]
The absolute value is bounded by $q\|n-1\|_2\|w\|_2$, so part~\textup{(ii)} of the lemma yields
\begin{equation}\label{eq:affine-mass}
 \big|\tfrac{d}{dt}\|w(t)\|_2\big|\le q\|n(t)-1\|_2\quad\text{a.e.}
\end{equation}
Since $w_0$ tends to zero at both ends of the line, $\|u_0\|_\infty=\|1+w_0\|_\infty\ge1$. Hence $q\int_0^{|t|}e^{q\|u_0\|_\infty^2s}\,ds\le M-1$. Using $\delta(t)\le\sqrt{8/q}\,\|\partial_xu_0\|_2M$ from \eqref{eq:cubic-density-explicit} and integrating \eqref{eq:affine-mass} yields
\begin{equation}\label{eq:affine-growth}
 \|w(t)\|_2\le\|w_0\|_2+
 \bigl(\|n_0-1\|_2+\sqrt{8/q}\,\|\partial_xu_0\|_2\bigr)(M-1).
\end{equation}
Also, \eqref{eq:cubic-current-growth} gives $\|w_x(t)\|_2\le\|\partial_xu_0\|_2+2\|\partial_xu_0\|_2(M-1)$. Combining the two component bounds by the triangle inequality in $\R^2$ gives
\[
 \|w(t)\|_{H^1}\le\|w_0\|_{H^1}
 +(M-1)\bigl[\|n_0-1\|_2+(2+\sqrt{8/q})\|\partial_xu_0\|_2\bigr],
\]
which proves \eqref{eq:cubic-affine-growth} and \eqref{eq:cubic-metric-growth} on affine balls.
\end{proof}

\section{Generalizing the net mass transfer method}\label{sec:current}

In Section~\ref{sec:cubic} we developed the net mass transfer argument for the cubic nonlinearity. Lemmas~\ref{lem:current} and~\ref{lem:relative} extend the net mass transfer estimates and the kinetic balance to other nonlinearities. We use these lemmas in the proofs of Propositions~\ref{prop:bounded-slope}--\ref{prop:subcubic-powers}.

\begin{lemma}[Generalized net mass transfer]\label{lem:current}
Let $I$ be an interval containing $0$, and let $u\in C(I;X^1)$ be a mild solution of
\[
iu_t+u_{xx}+Hu=0
\]
with $H\in C(I;X^1)$ real and $u(0)=u_0$. Denote
\[
 n=|u|^2,\quad n_0=|u_0|^2,\quad
 j=2\operatorname{Im}(\bar u u_x),\quad
 F(t)=-\int_0^t j(s)\,ds .
\]

\smallskip\noindent\textup{(i) Identities and amplitude estimate.}
Then $F\in C(I;H^2)\cap C^1(I;L^2)$, $F_x=n-n_0\in C(I;H^1)$ and
\begin{equation}\label{eq:weight}
 \int nF^2\,dx=\int n_0F^2\,dx\le \|u_0\|_\infty^2\|F\|_2^2 .
\end{equation}
The net mass transfer satisfies
\begin{equation}\label{eq:general-current-estimate}
 \left|\frac{d}{dt}\|F(t)\|_2\right|
 \le 2\|u_0\|_\infty\|u_x(t)\|_2
 \qquad\text{for almost every }t\in I.
\end{equation}
The kinetic energy and the squared density increments  satisfy
\begin{equation}\label{eq:current-balances}
 \frac{d}{dt}\|u_x(t)\|_2^2=\int H_xj\,dx,\qquad
 \tfrac12\frac{d}{dt}\|F_x\|_2^2=\int(n_x-\partial_x n_0)j\,dx .
\end{equation}
Moreover, for any fixed smooth nonnegative averaging kernel $\chi\in C_c^\infty((-1,1))$ with $\int\chi=1$,
\begin{equation}\label{eq:current-sup}
 \|u(t)\|_\infty\le \bigl(\|u_0\|_\infty^2+\|\chi'\|_2\|F(t)\|_2\bigr)^{1/2}+\|u_x(t)\|_2.
\end{equation}

\smallskip\noindent\textup{(ii) Conditional growth bounds.}
If, in addition, there are $L\ge0$ and a real measurable function $g(t,x)$ with $|g|\le L/2$ such that, for almost every $t\in I$,
\[
 Y=u_x+gFu,\qquad \|Y(t)\|_2\le \|\partial_xu_0\|_2+\tfrac L2\|u_0\|_\infty \|F(t)\|_2,
\]
then throughout $I$,
\begin{align}
 \|F(t)\|_2&\le 2\|u_0\|_\infty \|\partial_xu_0\|_2\int_0^{|t|}e^{L\|u_0\|_\infty^2s}\,ds,\label{eq:transfer-L}\\
 \|u_x(t)\|_2&\le \|\partial_xu_0\|_2\bigl(2e^{L\|u_0\|_\infty^2|t|}-1\bigr).\label{eq:main-L}
\end{align}
Together with \eqref{eq:current-sup}, these estimates bound $\|u(t)\|_{X^1}$ in terms of $\|u_0\|_\infty,\|\partial_xu_0\|_2,L$ and $|t|$.
\end{lemma}
\begin{proof}
\textbf{(i) Identities and amplitude estimate.}
Apply the spatial mollification argument in the proof of Lemma~\ref{lem:restricted-flows}(i) with forcing $Hu$. Since $H$ is real, $\operatorname{Im}(\bar u\,Hu)=0$, so $n_t=-j_x$ in distributions.

Now $j\in C(I;L^2)$, so $F_t=-j$ and $F_x=n-n_0$. We already have $F\in C(I;L^2)$, and $F_{xx}=n_x-\partial_x n_0\in C(I;L^2)$. Hence \eqref{eq:line-interpolation} gives $F\in C(I;H^2)$ and therefore $F^3\in W^{1,1}(\R).$ Now integration of $F_xF^2=(F^3)_x/3$ gives the mass-weighted cancellation \eqref{eq:weight}.

Testing $F_t=-j$ against $F$ and using \eqref{eq:weight} gives
\[
 \tfrac12\left|\frac{d}{dt}\|F(t)\|_2^2\right|
 =\left|\int Fj\dd\right|
 \le 2\|Fu\|_2\|u_x(t)\|_2
 \le 2\|F(t)\|_2\|u_0\|_\infty\|u_x(t)\|_2;
\]
now Lemma~\ref{lem:norm-differentiation}(ii) yields \eqref{eq:general-current-estimate}.

For the first identity in \eqref{eq:current-balances}, apply Lemma~\ref{lem:norm-differentiation}(i) to the spatially differentiated Duhamel formula with $v=u_x$ and $G=i(Hu)_x$. 
The second identity follows in the same way as Claim~\ref{cl:gradient-norm}, using $F_{xx}=n_x-\partial_xn_0$ and $F_t=-j$.

To prove \eqref{eq:current-sup}, fix $t\in I$ and $x\in\R$. Whenever $\chi(y-x)\ne0$, we have $|x-y|<1$, so the fundamental theorem of calculus and Cauchy--Schwarz give
\[
 |u(t,x)-u(t,y)|
 =\left|\int_y^x u_x(t,z)\,dz\right|
% \le |x-y|^{1/2}\|u_x(t)\|_2
 \le \|u_x(t)\|_2 \implies
 |u(t,x)|\le |u(t,y)|+\|u_x(t)\|_2.
\]
Since $\chi\ge0$ and $\int_{\R}\chi(y-x)\,dy=1$, averaging in $y$ yields
\[
 \begin{aligned}
 |u(t,x)|
 &=\int_{\R}\chi(y-x)|u(t,x)|\,dy\\
 &\le\int_{\R}\chi(y-x)\bigl(|u(t,y)|+\|u_x(t)\|_2\bigr)\,dy\\
 &=\int_{\R}\chi(y-x)|u(t,y)|\,dy+\|u_x(t)\|_2 .
 \end{aligned}
\]
Splitting the weight $\chi=\sqrt{\chi}\sqrt{\chi}$ and applying Cauchy--Schwarz, we obtain
\[
 \begin{aligned}
 \int_{\R}\chi(y-x)|u(t,y)|\,dy
 &\le\left(\int_{\R}\chi(y-x)\,dy\right)^{1/2}
      \left(\int_{\R}\chi(y-x)|u(t,y)|^2\,dy\right)^{1/2}\\
 &=\left(\int_{\R}\chi(y-x)|u(t,y)|^2\,dy\right)^{1/2}.
 \end{aligned}
\]
To bound this last integral, use $|u(t)|^2=|u_0|^2+F_x(t)$ and integrate by parts against the compactly supported kernel to obtain
\[
 \begin{aligned}
 \int_{\R}\chi(y-x)|u(t,y)|^2\,dy
 &=\int_{\R}\chi(y-x)|u_0(y)|^2\,dy-\int_{\R}\chi'(y-x)F(t,y)\,dy\\
 &\le \|u_0\|_\infty^2\int_{\R}\chi(y-x)\,dy
       +\|\chi'(\cdot-x)\|_2\|F(t)\|_2\\
 &=\|u_0\|_\infty^2+\|\chi'\|_2\|F(t)\|_2 .
 \end{aligned}
\]
Combining the estimates and taking the supremum over $x$ gives \eqref{eq:current-sup}.

\medskip\noindent\textbf{(ii) Conditional growth bounds.}
Reality of $gF$ gives $j=2\operatorname{Im}(\bar u Y)$. By the mass-weighted cancellation \eqref{eq:weight},
\[
 \tfrac12 \left|\frac{d}{dt}\|F(t)\|_2^2\right|
 =\left|\int Fj\,dx\right|
 \le 2\|Fu\|_2\|Y\|_2
 \le 2\|u_0\|_\infty \|F(t)\|_2\|Y\|_2 .
\]
The assumed bound on $\|Y\|_2$ and Lemma~\ref{lem:norm-differentiation}(ii) give
\[
 \left|\frac{d}{dt}\|F(t)\|_2\right|
 \le 2\|u_0\|_\infty \|\partial_xu_0\|_2+L\|u_0\|_\infty^2 \|F(t)\|_2
 \quad\text{a.e.}
\]
Since $F(0)=0$, Gronwall's inequality in both time directions gives \eqref{eq:transfer-L}. The triangle inequality and \eqref{eq:weight} then give
\[
 \begin{aligned}
 \|u_x(t)\|_2&\le\|Y\|_2+\tfrac L2\|u_0\|_\infty \|F(t)\|_2\\
 &\le \|\partial_xu_0\|_2+L\|u_0\|_\infty \|F(t)\|_2.
 \end{aligned}
\]
Substituting \eqref{eq:transfer-L} into this inequality gives \eqref{eq:main-L}.
\end{proof}

We now specialize to $H=f(|u|^2)-f(1)$. The kinetic energy can be compared with its initial value even when the potential energy is infinite.

\begin{lemma}[Relative kinetic balance]\label{lem:relative}
Let $f\in C^2([0,\infty);\R)$ and let $u\in C(I;X^1)$ be a mild solution of \eqref{eq:nls} on an interval containing zero. Write $n=|u|^2$, $n_0=|u_0|^2$, $j=2\operatorname{Im}(\bar u u_x)$ and $F(t)=-\int_0^t j(s)\,ds$.
For $0\le\theta\le1$, define the \emph{interpolated density} $n_\theta$ and the \emph{averaged density response} $\bar f$ by
\begin{equation}\label{eq:fvcdfgbhn}
 n_\theta=(1-\theta)n_0+\theta n,\qquad
 \bar f(n,n_0)=\int_0^1 f(n_\theta)\,d\theta,
\end{equation}
and, for any anti-derivative $V$ of $f$, put
\[
 P(n,n_0)=V(n)-V(n_0)-f(n_0)(n-n_0).
\]
Then
\begin{equation}\label{eq:f-balance}
 \|u_x(t)\|_2^2-\|\partial_xu_0\|_2^2
 =\int P(n,n_0)\dd-\int F(f(n_0))_x\dd.
\end{equation}
Equivalently, in terms of the averaged density response,
\begin{equation}\label{eq:averaged-kinetic-balance}
 \|u_x(t)\|_2^2-\|\partial_xu_0\|_2^2
 =-\int F\,\partial_x\bar f(n,n_0)\dd.
\end{equation}
\end{lemma}
\begin{proof}
Fix a compact interval $J\subset I$ containing $0 \in J$. Since $u\in C(I;X^1)$, we may set
\[
 M_J=\max_{t\in J}\|u(t)\|_\infty^2,\qquad L_J=\max_{0\le s\le M_J}|f'(s)|.
\]
Both $n(t,x)$ and $n_0(x)$ lie in $[0,M_J]$ for $t\in J$. Taylor's theorem with $V'=f$ gives
\[
 |P(n,n_0)|\le\tfrac12 L_J|n-n_0|^2 .
\]
Lemma~\ref{lem:current}(i) gives $F_x=n-n_0\in C(J;H^1)$, $F_t=-j$ and
\begin{equation}\label{eq:response-kinetic-derivative}
 \frac{d}{dt}\|u_x(t)\|_2^2=\int(f(n))_xj\dd .
\end{equation}
\begin{claim*}
The map $t\mapsto\int P(n(t),n_0)\dd$ is continuously differentiable on $J$, with
\begin{equation}\label{eq:relative-potential-derivative}
 \frac{d}{dt}\int P(n,n_0)\dd
 =\langle n_t,f(n)-f(n_0)\rangle_{H^{-1},H^1}
 =\int(f(n)-f(n_0))_xj\dd.
\end{equation}
\end{claim*}

\noindent\emph{Proof of the claim.} We first expand $P$ in its density variable. For $a,b,a+\delta\in[0,M_J]$, the Lipschitz bound on $f$ gives
\[
 \left|P(a+\delta,b)-P(a,b)-(f(a)-f(b))\delta\right|
 =\left|\int_0^\delta (f(a+s)-f(a))\,ds\right|
 \le \tfrac12 L_J|\delta|^2 .
\]
Fix $t$ and write $\Delta_h n=n(t+h)-n(t)$ for $t,t+h\in J$. Since $n-n_0=F_x$ and $F_t=-j\in C(J;L^2)$ by Lemma~\ref{lem:current}(i), $\Delta_h n/h\to -j_x(t)=n_t(t)$ in $H^{-1}$. Also $\Delta_h n\to0$ in $H^1$. Applying the density expansion pointwise with $a=n(t)$, $b=n_0$ and $\delta=\Delta_h n$, we obtain
\[
 \begin{aligned}
 &\left|\frac1h\int\bigl[P(n(t+h),n_0)-P(n(t),n_0)\bigr]\dd
 -\left\langle\frac{\Delta_h n}{h},f(n(t))-f(n_0)\right\rangle_{H^{-1},H^1}\right|\\
 &\qquad\le \frac{L_J}{2|h|}\|\Delta_h n\|_2^2
 \le \frac{L_J}{2}\left\|\frac{\Delta_h n}{h}\right\|_{H^{-1}}\|\Delta_h n\|_{H^1}
 \longrightarrow 0 .
 \end{aligned}
\]
By Lemma~\ref{lem:density-composition}, $f(n)-f(n_0)\in C(J;H^1)$. Thus the pairing converges to the one in \eqref{eq:relative-potential-derivative}, which is continuous in $t$. Substituting $n_t=-j_x$ and integrating by parts in space gives the last equality. This proves the claim.

Subtracting \eqref{eq:relative-potential-derivative} from \eqref{eq:response-kinetic-derivative} and integrating from zero to $t$ proves \eqref{eq:f-balance}.
To justify the integration by parts in \eqref{eq:averaged-kinetic-balance}, observe that the $n_\theta$ are uniformly bounded in $X^1$ for $t\in J$ and $0\le\theta\le1$. Since $n_\theta-n_\eta=(\theta-\eta)F_x\in H^1$, Lemma~\ref{lem:density-composition} gives
\[
 \|f(n_\theta)-f(n_\eta)\|_{H^1}
 \le C_J|\theta-\eta|\,\|F_x\|_{H^1},\qquad \theta,\eta\in[0,1].
\]
Taking $\eta=0$ gives $f(n_\theta)-f(n_0)\in H^1$, and the estimate for general $\theta,\eta$ gives Lipschitz continuity in $\theta$ with values in $H^1$. Its integral over $\theta\in[0,1]$ is thus $\bar f(n,n_0)-f(n_0)\in H^1$. Since $F\in H^1$ as well, integration by parts gives
\[
 \begin{aligned}
 &\int P(n,n_0)\dd-\int F(f(n_0))_x\dd\\
 &\quad=\int F_x\bigl(\bar f(n,n_0)-f(n_0)\bigr)\dd-\int F(f(n_0))_x\dd\\
 &\quad=-\int F\bigl(\partial_x\bar f(n,n_0)-(f(n_0))_x\bigr)\dd-\int F(f(n_0))_x\dd\\
 &\quad=-\int F\,\partial_x\bar f(n,n_0)\dd.
 \end{aligned}
\]
\end{proof}

\section{Proofs of Propositions~\ref{prop:bounded-slope}--\ref{prop:subcubic-powers}}\label{sec:lip}

\begin{proof}[\normalfont\bfseries Proof of Proposition~\ref{prop:bounded-slope}]
Gallo's local theory gives a unique maximal mild solution $u\in C(I;X^1)$, with $I=(T_-,T_+)\ni0$, and the blow-up alternative in the $X^1$ norm \cite[Theorems~4.1--4.2 and Proposition~4.1]{gallo}. We use the averaged density response $\bar f$ of \eqref{eq:fvcdfgbhn} to construct the modified derivative $Y$ and obtain global continuation. We then prove Lipschitz dependence and restrict the flow to $\E_1$ and $1+H^1$.

\medskip\noindent\textbf{Step 1. A priori bound and global continuation.}
Write $n=|u|^2$, $n_0=|u_0|^2$ and $F(t)=-\int_0^t2\operatorname{Im}(\bar u u_x)(s)\,ds$. Define the coefficients associated with the initial and current densities, respectively, by
\begin{equation}\label{eq:averaged-response-coefficients}
 g_1(n,n_0)=\int_0^1(1-\theta)f'(n_\theta)\,d\theta,\qquad
 g_2(n,n_0)=\int_0^1\theta f'(n_\theta)\,d\theta.
\end{equation}
Since $\partial_x n_\theta=(1-\theta)\partial_x n_0+\theta n_x$, differentiating the averaged density response gives
\[
 \partial_x\bar f=g_1\partial_xn_0+g_2n_x,
 \qquad |g_1|,|g_2|\le L/2.
\]
To apply Lemma~\ref{lem:current}(ii), we seek a modified derivative whose $L^2$ norm is controlled by the initial data and $\|F(t)\|_2$. Substituting the formula for $\partial_x\bar f$ into the kinetic balance \eqref{eq:averaged-kinetic-balance} and collecting the terms involving the current spatial derivatives on the left gives
\begin{equation}\label{eq:kinetic-balance-before-square}
 \|u_x\|_2^2+\int g_2F n_x\dd
 =\|\partial_xu_0\|_2^2-\int g_1F\partial_xn_0\dd .
\end{equation}
The left-hand side is the beginning of a square. Indeed, $n_x=2\operatorname{Re}(\bar u u_x)$ and $g_2F$ is real, so
\[
 g_2F n_x=2\operatorname{Re}\bigl(\overline{u_x}\,g_2Fu\bigr).
\]
Thus we choose $Y=u_x+g_2Fu$. Adding the missing square $\|g_2Fu\|_2^2$ to both sides of \eqref{eq:kinetic-balance-before-square} yields
\[
 \|Y\|_2^2
 =\|\partial_xu_0\|_2^2-\int g_1F\partial_xn_0\dd+\int g_2^2nF^2\dd .
\]
For the last term, first use $|g_2|\le L/2$ and then the mass-weighted cancellation \eqref{eq:weight} to obtain
\[
 \int g_2^2nF^2\dd
 \le \tfrac{L^2}{4}\int nF^2\dd
 =\tfrac{L^2}{4}\int n_0F^2\dd
 \le\tfrac{L^2}{4}\|u_0\|_\infty^2\|F(t)\|_2^2.
\]
Since $|g_1|\le L/2$ and $|\partial_xn_0|\le2\|u_0\|_\infty|\partial_xu_0|$, Cauchy--Schwarz bounds $\big|\int g_1F\partial_xn_0\dd\big|$ by $L\|u_0\|_\infty\|F(t)\|_2\|\partial_xu_0\|_2$. Combining these estimates gives
\[
 \|Y\|_2^2
 \le\left(\|\partial_xu_0\|_2+\tfrac L2\|u_0\|_\infty\|F(t)\|_2\right)^2.
\]
Applying Lemma~\ref{lem:current}(ii) with $g(t,x)=g_2(n(t,x),n_0(x))$ gives \eqref{eq:transfer-L} and \eqref{eq:main-L} and bounds the full $X^1$ norm on $I\cap[-T,T]$ in terms of $L,\|u_0\|_\infty,\|\partial_xu_0\|_2$ and $T$. The blow-up alternative therefore gives $I=\R$.

\medskip\noindent\textbf{Step 2. Lipschitz dependence in $X^1$.}
Fix $R,T>0$, and let $u,v$ be the solutions with initial data $u_0,v_0$ satisfying
\[
 \max\{\|u_0\|_{X^1},\|v_0\|_{X^1}\}\le R.
\]
Step~1 provides a common bound $B=B(L,R,T)\ge1$ such that
\[
 \sup_{|t|\le T}\max\{\|u(t)\|_{X^1},\|v(t)\|_{X^1}\}\le B.
\]
Applying the product rule and the derivative bound in the proof of Lemma~\ref{lem:density-composition} to $n=|u|^2$ and $b=|v|^2$ gives, for every $|t|\le T$,
\[
 \begin{aligned}
 \|\mathcal N(u(t))-\mathcal N(v(t))\|_{X^1}
 &\le C_B\|u(t)-v(t)\|_{X^1},\\
 C_B&=C\left[L(1+B^2)+B^4\sup_{0\le s\le B^2}|f''(s)|\right],
 \end{aligned}
\]
with a universal $C$. The difference of the Duhamel formulas is
\[
 u(t)-v(t)=S(t)(u_0-v_0)
 +i\int_0^tS(t-s)\bigl[\mathcal N(u(s))-\mathcal N(v(s))\bigr]\,ds.
\]
Using \eqref{eq:free-group-bound} and Gronwall as in Step~5 of the proof of Theorem~\ref{thm:main}, we obtain
\[
 \|u-v\|_{C([-T,T];X^1)}
 \le K_Te^{K_TC_BT}\|u_0-v_0\|_{X^1}.
\]
The a priori growth bound uses only $L$; the Lipschitz estimate also uses the bound on $f''$ over the solution's density range.

Since $f$ is real and $C^2$, $\mathcal N(z)=(f(|z|^2)-f(1))z$ is Lipschitz on bounded disks, vanishes on the unit circle and satisfies $\bar z\mathcal N(z)\in\R$. Lemma~\ref{lem:restricted-flows} and Steps~1--2 give persistence and Lipschitz dependence in $\E_1$ and $1+H^1$.
\end{proof}

\begin{proof}[\normalfont\bfseries Proof of Proposition~\ref{prop:anti-derivative-growth}]
Fix an initial norm ball $\|u_0\|_{X^1}\le R$, $R>0$.
Let $u\in C(I;X^1)$ be the maximal mild solution given by the local $X^1$ theory \cite{gallo}.

\begin{claim*}
For the relative potential of Lemma~\ref{lem:relative} there is a constant $C_R$ such that
\begin{equation}\label{eq:relative-potential-upper}
 P(n,b)\le C_R(n-b)^2,\qquad n\ge0,\quad 0\le b\le R^2 .
\end{equation}
\end{claim*}

\noindent\emph{Proof of the claim.} Set $R'=2R^2+1$. Using $V'=f$ yields
\[
 0\le n\le R' \implies P(n,b)=\int_b^n \bigl(f(s)-f(b)\bigr)\,ds
   \le \tfrac12\sup_{[0,R']}|f'|\,(n-b)^2,
\]
while $V(n)\le C_V(1+n^2)$ yields
\[
\begin{aligned}
 n>R' \implies P(n,b) & \le C_V(1+n^2)+|V(b)|+|f(b)|(n-b)\\
 & \le C'_R n^2\le 4C'_R(n-b)^2.
 \end{aligned}
\]
where
\[
C'_R=2C_V+\max_{[0,R^2]}|V|+\max_{[0,R^2]}|f|
\]
and of course $0\le n-b\le n \le 2(n-b)$ in the case $n>R'$.
Taking
\[
C_R=\max\!\left\{\sup_{[0,R']}|f'|,\;4C'_R\right\}
\]
gives \eqref{eq:relative-potential-upper}. This completes the proof of the claim.

\smallskip

In what follows we work for $t$ in  a compact subinterval of $I$ containing zero, and denote
$n=|u|^2$ and $n_0=|u_0|^2$. Integrating by parts and applying the mass-weighted cancellation \eqref{eq:weight} gives
\begin{equation}\label{eq:density-increment-by-transfer}
 \begin{aligned}
 \|n-n_0\|_2^2&=-\int F(n_x-\partial_x n_0)\dd\\
 &\le 2\|u_0\|_\infty\|F(t)\|_2
       \bigl(\|u_x(t)\|_2+\|\partial_xu_0\|_2\bigr).
 \end{aligned}
\end{equation}
Using \eqref{eq:relative-potential-upper} and $|(f(n_0))_x|\le 2C_R|u_0||\partial_xu_0|$ in \eqref{eq:f-balance} gives
\[
 \|u_x(t)\|_2^2-\|\partial_xu_0\|_2^2
 \le C_R\|n-n_0\|_2^2+2C_R\|u_0\|_\infty\|F(t)\|_2\|\partial_xu_0\|_2 .
\]
By \eqref{eq:density-increment-by-transfer},
\[
 \begin{aligned}
 &\bigl(\|u_x(t)\|_2-\|\partial_xu_0\|_2\bigr)
  \bigl(\|u_x(t)\|_2+\|\partial_xu_0\|_2\bigr)\\
 &\qquad\le 2C_R\|u_0\|_\infty\|F(t)\|_2
  \bigl(\|u_x(t)\|_2+2\|\partial_xu_0\|_2\bigr)\\
 &\qquad\le 4C_R\|u_0\|_\infty\|F(t)\|_2
  \bigl(\|u_x(t)\|_2+\|\partial_xu_0\|_2\bigr).
 \end{aligned}
\]
Dividing by $\|u_x(t)\|_2+\|\partial_xu_0\|_2$ when it is positive gives
\[
 \|u_x(t)\|_2\le\|\partial_xu_0\|_2+4C_R\|u_0\|_\infty\|F(t)\|_2 .
\]
Substituting the preceding kinetic bound into the current estimate \eqref{eq:general-current-estimate} gives
\[
 \begin{aligned}
 \left|\frac{d}{dt}\|F(t)\|_2\right|
 &\le 2\|u_0\|_\infty\|u_x(t)\|_2\\
 &\le 2\|u_0\|_\infty\|\partial_xu_0\|_2
     +8C_R\|u_0\|_\infty^2\|F(t)\|_2\quad\text{a.e.}
 \end{aligned}
\]
Since $F(0)=0$, Gronwall's inequality in both time directions and \eqref{eq:current-sup} bound $\|u(t)\|_{X^1}$ uniformly for $t\in I\cap[-T,T]$ and $\|u_0\|_{X^1}\le R$. The blow-up alternative gives $I=\R$.

For fixed $R,T>0$, let $B\ge1$ be the preceding bound for $\|u(t)\|_{X^1}$ on $[-T,T]$ for all $\|u_0\|_{X^1}\le R$. The corresponding densities take values in $[0,B^2]$. Thus the nonlinear difference estimate in Step~2 of the proof of Proposition~\ref{prop:bounded-slope} applies with $L=\sup_{[0,B^2]}|f'|$ and $\sup_{[0,B^2]}|f''|$, both finite since $f\in C^2$. The same Duhamel estimate and Gronwall's inequality give Lipschitz dependence in $X^1$. Lemma~\ref{lem:restricted-flows} now gives persistence and Lipschitz dependence in $\E_1$ and $1+H^1$, as in Proposition~\ref{prop:bounded-slope}.
\end{proof}

\begin{proof}[\normalfont\bfseries Proof of Proposition~\ref{prop:subcubic-powers}]
Put $\alpha=(p-1)/2 \in [1/2,1)$. For each $0<\varepsilon\le1$, the regularized response $f_\varepsilon(s)=q(s+\varepsilon)^\alpha$ is $C^2$ on $[0,\infty)$. Gallo's local theory therefore gives a unique local mild solution $u_\varepsilon$ of \eqref{eq:nls} with initial datum $u_0$ \cite[Proposition~4.1 and Theorem~4.1]{gallo}. Since $f_\varepsilon$ has bounded slope, Proposition~\ref{prop:bounded-slope} extends $u_\varepsilon$ globally. We obtain a kinetic bound uniform in $\varepsilon$ in Step~1 and construct the solution of \eqref{eq:subcubic-power} as the limit of $u_\varepsilon$ in Step~2, including at $p=2$. Without loss of generality we assume $\|u_0\|_\infty>0$.

\medskip\noindent\textbf{Step 1. The weighted estimate.}
Denote $n=|u_\varepsilon|^2$, $n_0=|u_0|^2$ and let $F$ be the net mass transfer function for $u_\varepsilon$. Define $g_1,g_2$ by \eqref{eq:averaged-response-coefficients} for $f_\varepsilon$. The interpolated density in \eqref{eq:fvcdfgbhn} satisfies $n_\theta\ge\theta n$ and $n_\theta\ge(1-\theta)n_0$. For $n>0$ and $0<\theta\le1$, the first bound and $\alpha-1<0$ give
\[
 f_\varepsilon'(n_\theta)
 =q\alpha(n_\theta+\varepsilon)^{\alpha-1}
 \le q\alpha(\theta n)^{\alpha-1}.
\]
Multiplying by $\theta$ and integrating gives
\[
 0\le g_2=\int_0^1\theta f_\varepsilon'(n_\theta)\,d\theta
 \le q\alpha n^{\alpha-1}\int_0^1\theta^\alpha\,d\theta
 =\frac{q(p-1)}{p+1}n^{(p-3)/2}.
\]
Using $n_\theta\ge(1-\theta)n_0$ and integrating with weight $1-\theta$ gives the same bound for $g_1$ with $n_0$ in place of $n$ wherever $n_0>0$. Thus, denoting $c_{p,q}=q(p-1)/(p+1)$ for brevity, we have $|g_2u_\varepsilon|^2\le c_{p,q}^2 n^{p-2}$ wherever $n>0$ and $|g_2u_\varepsilon|=0$ if $n=0$. For $2<p<3$, H\"older's inequality and \eqref{eq:weight} give
\[
 \|g_2Fu_\varepsilon\|_2^2\le c_{p,q}^2 \int n^{p-2}F^2\dd
 \le c_{p,q}^2 \left(\int nF^2\dd\right)^{p-2}
      \left(\int F^2\dd\right)^{3-p}
 \le c_{p,q}^2\|u_0\|_\infty^{2(p-2)}\|F(t)\|_2^2.
\]
At $p=2$, the pointwise estimate gives $|g_2u_\varepsilon|=g_2\sqrt n\le q/3$ wherever $n>0$. Since $g_2u_\varepsilon=0$ wherever $n=0$, 
\[
\|g_2Fu_\varepsilon\|_2\le \tfrac q3\|F(t)\|_2.
\]
Moreover, since $|\partial_x n_0|\le2|u_0||\partial_xu_0|$,  for $2\le p<3$ we get
\[
 \left|\int g_1F\partial_xn_0\dd\right|
 \le 2c_{p,q}\|u_0\|_\infty^{p-2}\|F(t)\|_2\|\partial_xu_0\|_2.
\]
Expanding the square of $Y=\partial_xu_\varepsilon+g_2Fu_\varepsilon$ and applying \eqref{eq:averaged-kinetic-balance} with $f=f_\varepsilon$ gives
\[
 \|Y\|_2^2=\|\partial_xu_0\|_2^2-\int g_1F\partial_xn_0\dd+\|g_2Fu_\varepsilon\|_2^2
 \le \bigl(\|\partial_xu_0\|_2+c_{p,q}\|u_0\|_\infty^{p-2}\|F(t)\|_2\bigr)^2.
\]
Since $j=2\operatorname{Im}(\bar u_\varepsilon Y)$, the same use of $F_t=-j$ and \eqref{eq:weight} as in Lemma~\ref{lem:current}(ii) gives
\[
 \left|\frac{d}{dt}\|F(t)\|_2\right|
 \le 2\|u_0\|_\infty\|\partial_xu_0\|_2
     +2c_{p,q}\|u_0\|_\infty^{p-1}\|F(t)\|_2
 \quad\text{a.e.}
\]
Gronwall's inequality with $F(0)=0$, followed by
$\|\partial_xu_\varepsilon\|_2\le\|Y\|_2+\|g_2Fu_\varepsilon\|_2$, gives \eqref{eq:subcubic-kinetic} for $u_\varepsilon$. The bound on $F$ and \eqref{eq:current-sup} also bound $\|u_\varepsilon(t)\|_{X^1}$ uniformly for $|t|\le T$ and initial data in each $X^1$ ball.

\medskip\noindent\textbf{Step 2. Passage to the power law, $\varepsilon \to 0$.}
Denote
\[
 \mathcal N_\varepsilon(z)=(f_\varepsilon(|z|^2)-f_\varepsilon(1))z, \qquad
 \mathcal N(z)=q(|z|^{p-1}-1)z.
\]
We identify $z\in\mathbb C$ with $(\Re z,\Im z)\in\R^2$ and regard $\mathcal N,\mathcal N_\varepsilon$ as maps from $\R^2$ to $\R^2$. %We write $D\mathcal N_\varepsilon$ for the real Jacobian and $D^2\mathcal N_\varepsilon$ for the second differential. 
Direct differentiation gives the componentwise bound
\[
 \bigl|\partial_i\partial_j(\mathcal N_\varepsilon)_k(z)\bigr|
 \le C_{p,q}(|z|^2+\varepsilon)^{(p-2)/2},
 \qquad i,j,k\in\{1,2\}.
\]
Since $p\ge2$, these second derivatives are uniformly bounded on each disk. The maps and their first derivatives converge pointwise to $\mathcal N$ and to $D\mathcal N$, respectively, including at zero. Equicontinuity therefore gives $C^1$ convergence on compact disks. In particular, at $p=2$ the uniform Hessian bound makes the Jacobians $D\mathcal N_\varepsilon$ Lipschitz on each disk with a common constant, and their limit $D\mathcal N$ inherits this bound. This is the regularity used in the derivative difference estimate below. By the chain rule, $\partial_x\mathcal N_\varepsilon(u)=D\mathcal N_\varepsilon(u)u_x$, where the Jacobian $D$ acts on the vector $(\Re u_x,\Im u_x)$. For $R>0$, each $u$ with $\|u\|_{X^1}\le R$ satisfies $\|u\|_\infty\le R$ and $\|u_x\|_2\le R$. Applying the chain rule to the difference gives
\[
 \sup_{\|u\|_{X^1}\le R}\|\mathcal N_\varepsilon(u)-\mathcal N(u)\|_{X^1}
\le \sup_{|z|\le R}|\mathcal N_\varepsilon(z)-\mathcal N(z)|
\,\,+\,\,R\sup_{|z|\le R}\|D\mathcal N_\varepsilon(z)-D\mathcal N(z)\|
\]
which tends to $0$ as $\varepsilon\to0.$ To obtain a common Lipschitz bound on each ball, use
\[
 \partial_x(\mathcal N_\varepsilon(u)-\mathcal N_\varepsilon(v))
 =D\mathcal N_\varepsilon(u)(u_x-v_x)
  +(D\mathcal N_\varepsilon(u)-D\mathcal N_\varepsilon(v))v_x .
\]
On $[-T,T]$, Step~1 places the regularized solutions in a common $X^1$ ball $B$. Let $L$ be the common Lipschitz constant of the maps $\mathcal N_\varepsilon$ on $B$. Subtracting the two mild forms and using \eqref{eq:free-group-bound} gives, for $0\le t\le T$,
\[
 \|u_\varepsilon(t)-u_\eta(t)\|_{X^1}
\le K_T L\int_0^t \|u_\varepsilon(s)-u_\eta(s)\|_{X^1}\,ds\\
\,\, +\,\, K_T t\sup_{v\in B}
       \|\mathcal N_\varepsilon(v)-\mathcal N_\eta(v)\|_{X^1}.
\]
The supremum tends to zero as $\varepsilon,\eta\to0$ by uniform convergence on $B$. Gronwall's inequality and the same estimate backward in time show that $(u_\varepsilon)$ is Cauchy in $C([-T,T];X^1)$. The limit solves \eqref{eq:subcubic-power} and inherits \eqref{eq:subcubic-kinetic}. The same Duhamel estimate for $\mathcal N$ gives uniqueness and Lipschitz dependence on bounded sets.

The limiting nonlinearity $\mathcal N(z)=q(|z|^{p-1}-1)z$ vanishes on the unit circle and $\bar z\mathcal N(z)=q(|z|^{p-1}-1)|z|^2$ is real. On a disk $|z|,|w|\le B$ with $B\ge1$, the mean value theorem gives $\big||z|^{p-1}-|w|^{p-1}\big|\le(p-1)B^{p-2}|z-w|$, also at $p=2$ by the reverse triangle inequality. Splitting the difference of the products therefore gives $|\mathcal N(z)-\mathcal N(w)|\le q(1+pB^{p-1})|z-w|$. Step~1 gives a common $X^1$ bound on $[-T,T]$ for each initial ball, which passes to the limit by Step~2. Together with Lipschitz dependence in $X^1$ from Step~2, Lemma~\ref{lem:restricted-flows} now gives persistence and Lipschitz dependence in both $\E_1$ and $1+H^1$.
\end{proof}

\end{document}